\documentclass{amsart}
\usepackage{enumerate,amsmath,amssymb}
\usepackage[mathcal]{euscript}
\numberwithin{equation}{section}
\newtheorem{theorem}[equation]{Theorem}
\newtheorem{proposition}[equation]{Proposition}
\newtheorem{lemma}[equation]{Lemma}
\newtheorem{corollary}[equation]{Corollary}

\theoremstyle{definition}
\newtheorem{example}[equation]{Example}
\newtheorem{definition}[equation]{Definition}
\newtheorem{assumption}[equation]{Assumption}

\theoremstyle{remark}
\newtheorem{remark}[equation]{Remark}
\def\C{\mathbb C}
\def\N{\mathbb N}
\def\R{\mathbb R}
\def\Z{\mathbb Z}
\def\dbar{d\hspace*{-0.08em}\bar{}\hspace*{0.1em}}
\def\embed{\hookrightarrow}
\def\eps{\varepsilon}
\def\interior#1{\overset{\,{}_{\circ}}{#1}}
\def\mass{\mathfrak m}
\def\st{\,|\,}
\def\ff{\mathit{ff}}
\def\ffi{\mathit{fi}}
\def\lb{\mathit{lb}}
\def\rb{\mathit{rb}}
\def\Dom{\mathcal D}
\def\A{\mathcal A}
\def\D{\partial}
\def\E{\mathcal E}
\def\F{\mathcal F}
\def\M{\mathcal M}
\def\O{\mathcal O}
\def\U{\mathcal U}
\def\T{\mathcal T}

\def\sym{\,{}^{b}\!\sigma}
\DeclareMathOperator{\Diff}{Diff}
\DeclareMathOperator{\Ind}{ind}
\DeclareMathOperator{\Tr}{Tr}
\DeclareMathOperator{\spec}{spec}
\def\OpSpaceA#1{\Psi^{#1}_{c}(M;\Lambda)}
\def\OpSpaceB#1#2{\Psi^{#1}_{#2}(M;\Lambda)}

\begin{document}

\title[Resolvents of cone pseudodifferential operators]%
{Resolvents of cone pseudodifferential operators,
asymptotic expansions and applications}
\author{Juan B. Gil}
\address{Penn State Altoona, 3000 Ivyside Park, Altoona, PA 16601-3760}
\email{jgil@psu.edu}
\author{Paul A. Loya}
\address{Department of Mathematics\\ Binghamton University\\
Binghamton, NY 13902}
\email{paul@math.binghamton.edu}
\maketitle
\begin{abstract}
 We study the structure and asymptotic behavior of the resolvent of
 elliptic cone pseudodifferential operators acting on weighted Sobolev 
 spaces over a compact manifold with boundary. We obtain an asymptotic 
 expansion of the resolvent as the spectral parameter tends to infinity, 
 and use it to derive corresponding heat trace and zeta function 
 expansions as well as an analytic index formula. 
\end{abstract}
\section{Introduction}
\label{Intro}

In this paper we study the structure and asymptotic behavior of the 
resolvent of elliptic cone \emph{pseudodifferential} operators acting on 
weighted Sobolev spaces over a compact manifold with boundary. Our results 
complete (and contain) the existing descriptions of the resolvent of a 
cone differential operator (on Sobolev spaces), and provide a first account 
on the structure of resolvents, heat kernels, and complex powers of 
pseudodifferential operators on manifolds with conic singularities. 

Resolvent and heat kernel asymptotics on conic manifolds have been 
studied by many authors since the seminal paper by Cheeger \cite{Ch79}. 
For certain classes of first and second order symmetric operators there are
contributions by Callias \cite{Ca83}, Callias and Uhlmann \cite{CaUh84},
Br\"uning and Seeley \cite{BrSe87, BrSe91}, and Mooers \cite{Mo99}, 
to mention just a few. In \cite{Le97}, Lesch generalized the techniques 
of Br\"uning and Seeley and obtained more general results for 
selfadjoint differential operators of arbitrary order.

Following Schulze's approach for the study of operators on manifolds with 
edges, see e.g. \cite{Sz89}, the first author developed a parameter-dependent 
calculus (cf. \cite{GiHeat01}) that describes the resolvent of an elliptic 
cone differential operator that is not necessarily selfadjoint. In
particular, he introduced the appropriate notion of parameter-dependent
ellipticity that guarantees the existence of the resolvent and provides
good norm estimates. 
In Section~\ref{Parametrix} we will show that this ellipticity 
condition is not only sufficient but also necessary.
Later in \cite{LoRes01,LoRes201}, following Melrose's approach \cite{RBM2}, 
the second author studied the resolvent of elliptic cone differential 
operators from a more geometric viewpoint. To this end, he developed a 
parameter-dependent calculus that gives a precise description of the 
Schwarz kernel of the resolvent, providing a more convenient framework 
to analyze heat kernels, zeta functions, and other geometric invariants, 
see e.g. \cite{GiLo01}. 

In the setting of resolvents of close extensions of a cone differential 
operator, there are recent results by Schrohe and Seiler \cite{SchSe03}, 
by Falomir, Muschietti, Pisani, and Seeley \cite{FMPS}, and by
Falomir, Muschietti, and Pisani \cite{FMP}.
More recently, Gil, Krainer, and Mendoza \cite{GKM2} proved the existence 
of the resolvent and sectors of minimal growth for the closed extensions 
of a general cone differential operator.  To the best of our knowledge, 
resolvents of elliptic cone pseudodifferential operators have not been 
studied before in any setting.

In this work we consider a cone pseudodifferential operator 
$A\in x^{-\mu}\Psi_b^{\mu}(M)$, where $M$ is a compact manifold with 
boundary, $x$ is a boundary defining function for $\partial M$, $\mu$ is 
a positive real number, and $\Psi_b^{\mu}(M)$ is the class of 
$b$-pseudodifferential operators of order $\mu$, as introduced by Melrose.
Our main goal is to give a precise description of the resolvent 
$(A-\lambda)^{-1}$ when $A$ satisfies the aforementioned 
parameter-dependent ellipticity on a sector $\Lambda\subset\C$. 
We obtain an asymptotic expansion in $\lambda$ as $|\lambda|\to\infty$, 
and use it to derive heat trace asymptotics and zeta function expansions.  
For this purpose, we extend the existing pseudodifferential calculi
introduced in \cite{LoRes01,LoRes201} and define two new classes of 
operators arising in the parametrix construction used to analyze the 
resolvent. 

As in the case of a differential operator, the construction of a good 
parameter-dependent parametrix of $A-\lambda$ is crucial to describe the 
fine structure of the resolvent and its asymptotic behavior in $\lambda$.
However, when the given operator is not differential but rather a
genuine pseudodifferential operator, for instance $\sqrt{\Delta}$,
the parametrix construction requires a more delicate analytic treatment.  
The general idea is to design a parameter-dependent pseudodifferential
calculus tailoring the new features of the operators into the geometry 
of their Schwartz kernels.

To illustrate the main technical difficulty in the parametrix construction
for the operator family $A-\lambda$, let us discuss the related (but much 
simpler) situation of an operator in the $b$-calculus. Given a 
parameter-elliptic $b$-differential operator $A$, one can construct a 
parametrix $B(\lambda)$ of $A-\lambda$ such that 
\begin{equation} \label{AB}
(A - \lambda) B(\lambda) = 1 + R(\lambda),
\end{equation}
where $R(\lambda)$ is in the calculus with bounds, of order $-\infty$, 
vanishing to infinite order as $|\lambda|\to\infty$ in $\Lambda$. 
For a $b$-pseudodifferential operator, the error term $R(\lambda)$ in 
\eqref{AB} can only be made to vanish to order $-1$ in the calculus with 
bounds. Nonetheless, this decay already implies that $R(\lambda) \to 0$ as 
$|\lambda| \to \infty$, thus $1 + R(\lambda)$ can be inverted for large 
$\lambda$, and consequently, the resolvent exists and belongs to the calculus. 

However, when $A$ is a cone pseudodifferential operator, the additional 
weight factor $x^{-\mu}$ in $A$ makes the situation more complicated: 
One can obtain an expression similar to \eqref{AB}, but the boundary 
defining function $x$, the spectral parameter $\lambda$, and 
the bounds, are all coupled in a way that the operator $1 + R(\lambda)$ is 
unfortunately \emph{not} invertible even for large $\lambda$. 
A novelty of this paper is the development of two new parameter-dependent
calculi with bounds which incorporate the coupling of the boundary defining 
function, the spectral parameter, and the bounds. We introduce these
operator classes and show the corresponding composition theorems. 
This will allow us to further modify $R(\lambda)$ and get a true residual 
term that decays as $|\lambda|\to \infty$, so that $1 + R(\lambda)$ can be 
inverted within the calculus.

Once the resolvent of an elliptic cone pseudodifferential operator is
understood, we use its structure to study the corresponding heat kernels
and complex powers. In particular, the short-time asymptotic expansion of 
the heat trace obtained in this paper is used to get part of an analytic 
index formula consisting of two terms; a term coming from the heat trace 
asymptotics of an associated operator with no boundary spectrum, and a 
second term that resembles the eta invariant. This formula relies on an 
index formula by Piazza~\cite{Pi93} and on a factorization theorem 
proposed by Schulze and proved by Witt~\cite{Wi01}.   

We now outline the content of this paper. We begin in Section~\ref{sec:mwc} 
by reviewing various conormal spaces of functions on manifolds with
corners as introduced in Melrose's seminal paper \cite{MeR92}. With
this background, in Section~\ref{PsdoCalculus} we define and discuss the
new parameter-dependent pseudodifferential calculi that are needed
in Section~\ref{Parametrix} to construct a good parametrix for a 
parameter-elliptic cone pseudodifferential operator. 
In Section~\ref{AsymptoticExpansion} we use the structure of these 
calculi to obtain resolvent, heat kernel, and zeta function expansions. 
Finally, in Section~\ref{IndexFormula}, we discuss the index of the 
closure of an elliptic cone operator.

\medskip
{\bf Acknowledgments.} 
We thank Gerardo Mendoza for many valuable discussions, especially for 
those concerning the last section of this paper.

\section{Manifolds with corners, asymptotics, and $b$-operators}
\label{sec:mwc}

An $n$-dimensional manifold with corners $Z$ is a topological
space with $C^\infty$ structure given by local charts of the form
$[0,1)^k \times (-1,1)^{n-k}$, where $k$ can run between $0$ and
$n$ depending on where the chart is located in the manifold.
Each boundary hypersurface $H$ is embedded and has a globally defined
boundary defining function; a nonnegative function in $C^\infty(Z)$
that vanishes only on $H$ where it has a nonzero differential.

\subsection*{Asymptotic expansions}
Let $\U = [0,1)^k_x \times (-1,1)^{n-k}_y$. Then for $a \in \mathbb{R}^k$
the space of symbols $\Sigma^a(\U)$ consists of those smooth functions
$u\in C^\infty(\interior{\U})$ of the form
\[ u(x,y) = x_1^{a_1}\cdots x_k^{a_k}\, v(x,y), \]
where for each $\alpha$ and $\beta$, $(x\D_x)^\alpha \D_y^\beta
v(x,y)$ is a bounded function.

Let $\N$ be the set of positive integers and let $\N_0=\N\cup\{0\}$.
An \emph{index set} $E$ is a discrete subset of $\C \times \N_0$ such that
\begin{itemize}
\item $(z,k) \in E \Rightarrow (z,\ell) \in E$ for all
      $0\leq \ell\leq k$, and
\item given any $N \in \R$, the set $\{(z,k)\in E \st \Re z \leq N\}$
      is finite.
\end{itemize}
If in addition, $(z,k) \in E \Rightarrow (z+\ell,k) \in E$ for all
$\ell \in \N_0$, then $E$ is called a $C^\infty$ index set. For simplicity,
we will use the words ``index set'' instead of ``$C^\infty$ index set'' unless
stated otherwise. A discrete subset $D\subset\C$ will be referred to as
an index set by means of the identification $D \cong \{(z,0) \st z\in D\}$.

Given an index set $E$, a function $u\in \Sigma^a(\U)$ is said to have an
\emph{asymptotic expansion} at $x_1=0$ with index set $E$ if, for each $N>0$,
\begin{equation} \label{expprop}
 u(x,y) = \hspace{-1em} \sum_{(z,k) \in E,\, \Re z \leq N} \hspace{-1ex}
 x^z_1\, (\log x_1)^k\, u_{(z,k)}(x',y)\ +\ x_1^{N}\,u_N(x,y)
\end{equation}
with $u_N(x,y)\in \Sigma^a(\U)$ and $u_{(z,k)}(x',y)\in \Sigma^{a'}(\U')$,
where $a=(a_1,a')$, $x = (x_1,x')$, and
$\U' = [0,1)^{k-1}_{x'} \times (-1,1)^{n-k}_y$.  Furthermore, given
$\kappa >0$, the function $u$ is said to have a \emph{partial
expansion} at $x_1 = 0$ with index set $E$ of order $\kappa$ if $u$
admits the expansion \eqref{expprop} for all $N \leq \kappa$. In fact,
it is sufficient to check that \eqref{expprop} holds for $N = \kappa$. 
Observe that a function has an asymptotic expansion at $x_1 = 0$ with index 
set $E$ if and only if it has a partial expansion at $x_1 = 0$ with index 
set $E$ of any order $\kappa>0$.  Note also that if $E = \varnothing$, then 
the expansion property \eqref{expprop} holds for $N=\kappa$ if and only if 
$u$ vanishes to order $\kappa$ at $x_1 = 0$.  Asymptotic and partial 
asymptotic expansions at any other boundary $x_i = 0$ are defined similarly.

On a manifold with corners $Z$ one can define asymptotic expansions at a
hypersurface $H$ with index set $E$ by reference to local coordinates.
First of all, a function $u\in C^\infty(\interior Z)$ is said to be in
$\Sigma^0(Z)$, if for any patch $\U$ on $Z$ and for any $\varphi \in
C_c^\infty(\U)$, the function $\varphi u$ is an element of $\Sigma^0(\U)$.
Let $H_1,\ldots, H_m$ be the hypersurfaces of $Z$ with
corresponding boundary defining functions $\rho_1,\ldots,\rho_m$.
For $a\in\R^m$ we define
\begin{equation*}
\Sigma^a(Z) = \{ \rho_1^{a_1}\cdots\rho_m^{a_m}\,v \st v\in\Sigma^0(Z)\}.
\end{equation*}
A function $u\in \Sigma^a(Z)$ has a partial expansion at $H$ with index set
$E$ of order $\kappa$, if for any patch $\U = [0,1)_{x_1} \times \U'$ on
$Z$ with $H \cap \U = \{x_1 = 0\}$, and for any $\varphi \in C_c^\infty(\U)$,
the function $\varphi u$  has a partial expansion at $x_1 = 0$
with index set $E$ of order $\kappa$ in the sense described above.

If $\E$ is a collection of index sets $\E = \{E_{H_1},\ldots,E_{H_\ell}\}$
corresponding to some family of hypersurfaces $H_1,\ldots,H_\ell$ of $Z$,
then we denote by $\A^\E_\kappa(Z)$ the space of functions
$u\in \Sigma^a(Z)$ for some $a\in\R^m$ such that for each $H$, $u$ has a
partial expansion at $H$ with index set $E_H$ of order $\kappa$.
Finally, we define
\begin{equation*}
  A^\E(Z) = \bigcap_{\kappa>0}\A^\E_\kappa(Z).
\end{equation*}

\subsection*{Blow-up and pseudodifferential operators}
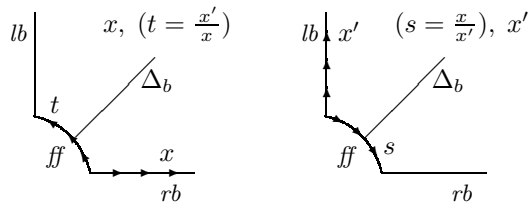
\begin{figure} \centering
\setlength{\unitlength}{.25mm} \centering
\begin{picture}(300,100)(-150,10) 
\put(-155,0){ \begin{picture}(120,120)(-30,-30) 
\put(-21,20){ \line(0,1){55}} \put(9,-10){ \line(1,0){55}} 
\qbezier(-17,20)(6,15)(13,-10)
\put(20,65){$x,\; (t= \frac{\,x'}{x})$} 
\put(9,0){\vector(-2,3){1}} \put(2,10){\vector(-1,1){1}} 
\put(-9,18){\vector(-3,2){1}}
\put(30,-10){\vector(1,0){1}} \put(45,-10){\vector(1,0){1}}
\put(60,-10){\vector(1,0){1}} 
\put(-8,22){$t$} \put(50,-4){$x$}
\put(50,-25){\small \emph{rb}} \put(-10,-5){\small \emph{ff}}
\put(-30,60){\small \emph{lb}} \put(1,9){ \line(1,1){42}}
\put(40,35){$\Delta_b$}
\end{picture} } 
\put(0,0){ \begin{picture}(120,120)(-30,-30)
\put(-21,20){ \line(0,1){55}} \put(9,-10){ \line(1,0){55}} 
\qbezier(-17,20)(6,15)(13,-10)
\put(20,65){$ (s = \frac{\!\! x}{x'}),\; x'$} 
\put(10,0){\vector(2,-3){1}} \put(2,10){\vector(1,-1){1}} 
\put(-8,17){\vector(2,-1){1}}
\put(-16,35){\vector(0,1){1}} \put(-16,50){\vector(0,1){1}}
\put(-16,65){\vector(0,1){1}} 
\put(-10,60){$x'$} \put(14,0){$s$}
\put(50,-25){\small \emph{rb}} \put(-10,-5){\small \emph{ff}}
\put(-30,60){\small \emph{lb}} \put(0,9){ \line(1,1){42}}
\put(40,35){$\Delta_b$}
\end{picture} } 
\end{picture}
\caption{Each of these coordinates together with coordinates on $Y^2$
define projective coordinates on $M^2_b$ near $\ff$.} \label{fig:proj}
\end{figure}

Let $M$ be a smooth manifold of dimension $n$ with boundary $Y =\partial M$. 
Then the product $M^2=M\times M$ is a manifold with corners in 
the above sense.  The blow-up $M^2_b=[M^2;Y^2]$ of $M^2$ along $Y^2$
(cf. \cite{RBM2}) is then a new manifold with corners that has an atlas
consisting of the usual coordinate patches on $M^2 \setminus Y^2$
together with polar coordinate patches over $Y^2$ in $M^2$.
For instance, if $M^2 =[0,\infty)_x\times [0,\infty)_{x'}$, then $M^2_b$
is the set $[0,\infty)_r\times (\mathbb{S}^1\cap M^2)_\theta$ with
$(r,\theta)=(\|(x,x')\|, \tan^{-1}(x'/x))$.
In this paper we will work with the more convenient projective
coordinates $(x,x') \mapsto (x,t)$ with $t=x'/x$, or $(x,x') \mapsto
(s,x')$ with $s=x/x'$. The boundary hypersurfaces $\lb$,
$\rb$, and $\ff$ (for ``left boundary'', ``right boundary'',
and ``front face'', respectively) of $M^2_b$ together with the
projective coordinates are shown in Figure~\ref{fig:proj}.

Henceforth we fix a $b$-measure $\mass$ on $M$ and we
denote by $\mass'$ the lift of $\mass$ to $M^2$ under the right
projection $M^2 \ni (x,x') \mapsto x'\in M$. 

\begin{definition}\label{b-Operators}
For $\mu \in \R$, the space $\Psi^\mu_b(M)$ of $b$-pseudodifferential 
operators consists of operators $A$ on $C^\infty(M)$ that have a 
Schwartz kernel $K_{A}$ satisfying the following two conditions:
\begin{itemize}
\item
Given $\varphi \in C^\infty_c(M^2_b \setminus \Delta_b)$, the kernel
$\varphi K_{A}$ is of the form $k\cdot\mass'$, where $k$ is a smooth function
on $M^2_b$ that vanishes to infinite order at the boundaries $\lb$ and $\rb$.
\item
Given a coordinate patch of $M^2_b$ overlapping $\Delta_b$ of the form 
$\U_y\times \R^n_z$ such that $\Delta_b \cong \U \times \{0\}$, and given 
$\varphi\in C^\infty_c(\U \times \R^n)$, we have
\[ (\varphi K_{A}) = \int e^{iz \cdot \xi}\, a(y,\xi)\,\dbar \xi\cdot\mass', \]
where $y \mapsto a(y,\xi)$ is smooth with values in $S^{\mu}_{c\ell}(\R^n)$;
the space of classical symbols of order $\mu$.
\end{itemize}
\end{definition}

The space $\Diff_b^m(M)$ of totally characteristic differential operators
of order $m$ is clearly contained in $\Psi^m_b(M)$.

\section{Parametric pseudodifferential calculus}
\label{PsdoCalculus}

The spaces of parametric symbols and pseudodifferential operators
discussed in this section are intended to describe operator families
of the form $B (A - \lambda)^{-1}$, where $A$ and $B$ are both
cone pseudodifferential operators on an compact manifold $M$, and
$\lambda$ is a spectral parameter living on a sector $\Lambda\subset\C$.
Our symbol calculus is somewhat related to the weakly parametric calculus 
from Grubb and Seeley \cite{GrSe95}.

\subsection*{Symbol spaces}

For $\mu$, $p\in\R$ and $d>0$ we define $S^{\mu,p,d}(\R^n;\Lambda)$ as
the space of functions $a\in C^\infty(\R^n\times \Lambda)$ such that
\[ |\partial_\xi^\alpha \partial_\lambda^\beta a(\xi,\lambda)|
 \le C_{\alpha \beta} (1+|\xi|)^{\mu-p-|\alpha|}
 (1+|\xi|+|\lambda|^{1/d})^{p-d|\beta|}. \]
The space $S_{r}^{\mu,p,d}(\R^n;\Lambda)$, $p/d\in\mathbb Z$, consists
of elements $a\in S^{\mu,p,d}(\R^n;\Lambda)$ such that if we set
\[ \tilde a(\xi,z):= z^{p/d}a(\xi,1/z), \]
then $\tilde{a} (\xi,z)$ is smooth at $z=0$, and
\begin{equation} \label{estimater}
|\partial_\xi^\alpha \partial_z^\beta \tilde a(\xi,z)| \le
C_{\alpha \beta} (1+|\xi|)^{\mu-p-|\alpha|+d|\beta|}
(1+|z||\xi|^{d})^{p/d-|\beta|} \end{equation}
uniformly for
$|z|\le 1$. Further let $S_{r,c\ell}^{\mu,p,d}(\R^n;\Lambda)$ be
the space of elements $a\in S_r^{\mu,p,d}(\R^n;\Lambda)$ that, for
every $N\in\N$, admit a decomposition
\begin{equation}
a(\xi,\lambda)=\sum_{j=0}^{N-1}
\chi(\xi)a_{\mu-j}(\xi,\lambda) + r_{N}(\xi,\lambda) \label{symbexp},
\end{equation}
where $r_N\in S_r^{\mu-N,p,d}(\R^n;\Lambda)$, $\chi\in C^\infty(\R^n)$
with $\chi(\xi)=0$ for $|\xi|\le \frac12$ and $\chi(\xi)=1$ for
$|\xi|\ge 1$, and where each $a_{\mu-j}(\xi,\lambda)$ 
has the following properties:
\begin{itemize}
\item  $a_{\mu-j}(\delta\xi,\delta^d\lambda)=
\delta^{\mu-j}a_{\mu-j}(\xi,\lambda)$ for every $\delta>0$,
\item $z^{p/d} a_{\mu-j}(\xi,1/z)$ is smooth at $z=0$.
\end{itemize}

\begin{example} \label{ex}
Let $a(\xi)$ be a homogeneous function in $\xi\in\R^n$ of degree $\mu
\in \R$ that never takes values in a sector $\Lambda$ for $\xi \ne 0$,
and let $b(\xi)$ be a homogeneous function of degree $\mu' \in \R$.
Given $\ell \in \N_0$, set
\[ q(\xi,\lambda) = b(\xi) (a(\xi) - \lambda)^{-\ell}.\]
Then, $\chi(\xi)\, q(\lambda,\xi) \in
S_{r,c\ell}^{\mu'-\ell\mu,-\ell\mu,\mu}(\R^n;\Lambda)$. Here, the cut-off
function $\chi(\xi)$ is needed because $a(\xi)$ and $b(\xi)$ are,
in general, not smooth at $\xi = 0$.
\end{example}

\subsection*{Parameter-dependent operators}

We first review some spaces of parameter-dependent pseudodifferential
operators used to capture resolvents of cone {differential}
operators (see \cite{LoRes01} and \cite{LoRes201}). Henceforth we
shall fix a boundary defining function $\varrho$ for $\ff$. Recall
that $\mass'$ denotes the fixed $b$-measure $\mass$ lifted to $M^2$
under the right projection $M^2 \ni (x,x') \mapsto x' \in M$.

\begin{definition}\label{SpaceA}
Given $\mu,p,d\in \R$ with $p/d \in \Z$ and $d>0$, the space
$\OpSpaceA{\mu,p,d}$ consists of parameter-dependent
operators $A(\lambda)$ that have a Schwartz kernel
$K_{A(\lambda)}$ satisfying the following two conditions:
\begin{itemize}
\item Given $\varphi \in C^\infty_c(M^2_b \setminus \Delta_b)$,
the kernel $\varphi K_{A(\lambda)}$ is of the form $k(\varrho^d
\lambda,q)\cdot\mass'$, where $k(\lambda,q)$ is a smooth function of
$(\lambda,q) \in \Lambda \times M^2_b$ that vanishes to infinite
order in $q$ at the sets $\lb$ and $\rb$, and is
such that if we define $\widetilde{k}(z,q)= z^{p/d}k(1/z, q)$,
then $\widetilde{k}(z,q)$ is smooth at $z=0$.
\item Given a coordinate patch of $M^2_b$ overlapping $\Delta_b$ of
the form $\U_y \times \R^n_\zeta$ such that $\Delta_b \cong\U \times \{0\}$,
and given $\varphi \in C^\infty_c(\U \times \R^n)$, we have
\[ \varphi K_{A(\lambda)} = \int e^{i\zeta \cdot \xi}\, 
   a(y,\xi,\varrho^d\lambda)\, \dbar \xi\cdot \mass', \]
where $y \mapsto a(y,\xi,\lambda)$ is smooth with values in
$S^{\mu,p,d}_{r,c\ell}(\R^n;\Lambda)$.
\end{itemize}
\end{definition}

Let $[\Lambda;\{0\}]$ be the sector $\Lambda$ blown-up at at the
origin; that is, $\Lambda$ with polar coordinates taken at
$\lambda = 0$, let $\overline{\Lambda}$ denote the
stereographic compactification of $[\Lambda;\{0\}]$ in the Riemann
sphere. Coordinates on $\overline{\Lambda}$ near the blown-up
origin are $\rho_0= |\lambda|$ and $\theta = \lambda/|\lambda|$;
near $\lambda = \infty$ the coordinates are $\rho_\infty = |\lambda|^{-1}$
and $\theta = \lambda/|\lambda|$. Let $d > 0$ and let
$\overline{\Lambda}_d=\{\lambda^{1/d} \st \lambda\in \overline{\Lambda}\}$
so that the radial coordinates on $\overline{\Lambda}_d$ are
$r_0= |\lambda|^{1/d}$ near the origin and $r_\infty = |\lambda|^{-1/d}$
near infinity.

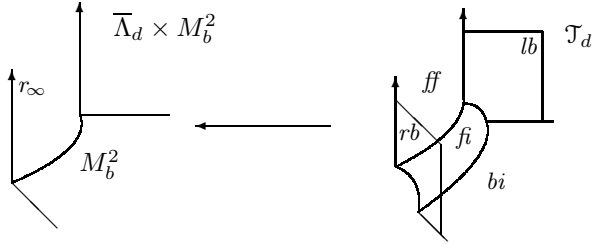
\begin{figure} \centering
\setlength{\unitlength}{.3mm} \centering 
\begin{picture}(340,90)(-170,-5)
\put(30,40){\vector(-1,0){60}}
\put(-120,20){ \begin{picture}(100,80)(-5,-25) 
\qbezier(0,-30)(35,-15)(30,0) \put(0,-30){\vector(0,1){50}}
\put(30,0){\vector(0,1){50}} 
\put(0,-30){\line(1,-1){20}} \put(30,0){\line(1,0){40}}
\put(3,10){\small $r_{\!\infty}$} \put(30,-25){$M^2_b$} 
\put(45,35){$\overline{\Lambda}_d \times M^2_b$}
\end{picture} }

\put(50,20){ \begin{picture}(100,80)(-5,-22) 
\qbezier(0,-20)(30,-5)(30,10) \qbezier(10,-40)(45,-15)(40,0)
\qbezier(0,-20)(12,-25)(10,-40) \qbezier(30,8)(38,8)(40,0)
\put(0,-20){\vector(0,1){40}} \put(30,10){\vector(0,1){40}}
\put(10,-40){\line(1,-1){13}} 
\put(40,0){\line(1,0){30}} \put(0,10){\line(1,-1){20}} 
\put(20,-50){\line(0,1){40}}
\put(30,40){\line(1,0){35}} \put(65,0){\line(0,1){40}}
\put(75,35){$\mathcal{T}_d$} \put(40,-30){\small \emph{bi}}
\put(27,-10){\small \emph{fi}} \put(12,15){\small \emph{ff}}
\put(1,-8){\small \emph{rb}} \put(55,30){\small \emph{lb}}
\end{picture} }
\end{picture}
\caption{The manifold $\T_d$ near infinity. Here, $r_{\infty} = 0$
defines the boundary at $|\lambda|=\infty$. }\label{fig:M2bL}
\end{figure}

We consider (see Figure \ref{fig:M2bL})
\begin{equation} \label{def:Td}
 \T_d:=[\overline{\Lambda}_d \times M^2_b;
 \{r_\infty=0\} \times \ff_{\!b}],
\end{equation}
the blow-up of $\overline{\Lambda}_d \times M^2_b$ along
$\{r_\infty=0\} \times \ff_{\!b}$, where
$\ff_{\!b}$ is the front face of $M^2_b$.

The blown-up manifold $\T_d$ has eight boundary hypersurfaces, five of
which are illustrated in Figure~\ref{fig:M2bL}, namely, $\ffi$ (face at
infinity), \emph{bi} (boundary at infinity), and the three hypersurfaces
$\lb$, $\rb$, and $\ff$, induced by the corresponding boundaries of the
manifold $M^2_b$. The other three hypersurfaces are $\{r_0=0\}$ and the
endpoints of the angular variable $\theta$. Because we are interested
in asymptotics for $|\lambda|$ near infinity, these three
hypersurfaces will play only minor roles.

\begin{assumption}
From now on, all functions depending on $\lambda$, either
implicitly (as functions on $\T_d$, for instance) or explicity (as
functions on $\overline{\Lambda}$), are assumed to be smooth in
$\theta=\lambda/|\lambda|$ and symbols of order zero at $\{r_0 = 0\}$.
\end{assumption}

We now use the notation from Section~\ref{sec:mwc} to describe the
various residual classes of pseudodifferential operators with asymptotics
that appear in the parametrix construction of parameter-dependent elliptic
operators.

\begin{definition}\label{AsympSpaceA}
Let $\E = (E_{\lb},E_{\rb},E_{\ff},E_{\ffi},\varnothing)$
be an index family for $\T_d$ associated to the faces
$(\lb,\rb,\ff,\ffi,bi)$. We denote by
$\OpSpaceA{-\infty,d,\E}$ the space of those parameter-dependent
operators $A(\lambda)$ that have a Schwartz kernel of the form
$K_{A(\lambda)} = k\cdot \mass'$ where $k \in \A^{\E}(\T_d)$.
Thus $k$ defines a function on $\T_d$ that vanishes to infinite order
at $bi$ and have asymptotic expansions at the hypersurfaces $\lb$,
$\rb$, $\ff$, and $\ffi$, determined by the index sets
$E_{lb}$, $E_{rb}$, $E_{\ff}$, and $E_{\ffi}$, respectively.
\end{definition}

\subsection*{Two new parameter-dependent residual classes}

In order to capture the resolvents of elliptic \emph{pseudodifferential}
operators we need to introduce two new classes of smoothing operators
satisfying only conormal bounds.  We begin by recalling the calculus
with bounds (cf. \cite[Section~5.16]{RBM2}).

\begin{definition}
Let $\alpha > 0$ and let $\E = (\varnothing, \varnothing, \N_0)$ be
the index family on $M^2_b$ associating the empty sets to its left and
right boundaries, and $\N_0$ to its front face. The space
$\Psi_b^{-\infty,\alpha}(M)$ denotes the class of operators $A$
having a Schwartz kernel of the form $K_A = k\cdot\mass'$, with 
$k\in \A^\E_{\alpha+\eps}(M^2_b)$ for some $\eps > 0$.
More precisely, if $\rho_{l}$ and $\rho_{r}$ are boundary defining
functions for the left and right boundaries of $M^2_b$, then the function
$\rho_{l}^{-\alpha - \eps}\!\rho_{r}^{-\alpha-\eps}\, k$
is a symbol in $\Sigma^0(M^2_b)$ having a partial expansion at the front
face of $M^2_b$ with index set $\N_0$ of order $\alpha+\eps$.
\end{definition}

\begin{definition}\label{SpaceB}
Let $N\in \N$ and $d>0$.
For $m\in\N$ we define $\OpSpaceB{-\infty,d}{m,N}$ as the space of
those parameter-dependent operators $A(\lambda)$ whose Schwartz kernel
$K_{A(\lambda)}$ is of the form $k(\varrho^d \lambda, q)\cdot\mass'$
with $k(\lambda,q)$ satisfying the following properties:
\begin{itemize}
\item[(a)] For some $\eps>0$,
 $\rho_{l}^{-Nd - \eps}\!\rho_{r}^{-Nd-\eps} k$ is a symbol in
 $\Sigma^0(\overline{\Lambda}\times M^2_b)$ having a partial expansion
 at the face $\overline{\Lambda}\times\ff$ with index set $\N_0$ of order
 $Nd + \eps$. Again, $\rho_{l}$ and $\rho_{r}$ are boundary defining
 functions for $\lb$ and $\rb$ in $M^2_b$,
\item[(b)] For each $N'\le N$,
\begin{equation*}
  k(\lambda,q) = \sum_{j=m}^{N' - 1} \lambda^{-j}\, f_j(q) +
  \lambda^{-N'} k_{N'}(\lambda,q),
\end{equation*}
where $f_j\in \A^{\E}_{2Nd-jd}(M^2_b)$ with
$\E = (\varnothing, \varnothing, \N_0)$,
and $k_{N'}$ satisfies (a) with $Nd$ replaced by $2Nd-N'd$. If $m\ge N$,
then we disregard the summation and require instead
$k(\lambda,q)=\lambda^{-N} k_{N}(\lambda,q)$, where $k_{N}$ satisfies (a).
\end{itemize}
\end{definition}

The next lemma relates the two parameter-dependent spaces introduced in
Definitions \ref{SpaceA} and \ref{SpaceB};
the proof follows immediately from the definitions.

\begin{lemma} \label{lem-regbounds}
If $p/d \in -\N$, then for any $N \in \N$,
\[ \OpSpaceA{-\infty,p,d} \subset \OpSpaceB{-\infty,d}{m,N},
   \quad m = - p/d. \]
\end{lemma}

Our next space of operators is a calculus with bounds version of the
space $\OpSpaceA{-\infty,d,\E}$ introduced in Definition~\ref{AsympSpaceA}.

\begin{definition}
Let $\E = (E_{\lb},E_{\rb},E_{\ff},E_{\ffi})$ be an index
family for $\T_d$ associated to the faces $(\lb,\rb,\ff,\ffi)$.
Then we define $\OpSpaceB{-\infty,d,\E}{N}$ as those parameter-dependent
operators $A(\lambda)$ that have a Schwartz kernel of the form
$K_{A(\lambda)} = k\cdot\mass'$, where $k$ is a symbol on $\T_d$, of order
$Nd$ at $bi$, that satisfies:
\begin{itemize}
\item Given $\varphi\in C^\infty(\T_d)$ supported near $\ffi$,
$\varphi\, k$ is in $\A^{\E}_{Nd + \eps}(\T_d)$ for some $\eps > 0$.
\item Given $\psi\in C^\infty(\T_d)$ supported away from
$\ffi$, $\psi\, k$ is the kernel of a parameter-dependent
operator in $\OpSpaceB{-\infty,d}{N,N}$.
\end{itemize}
\end{definition}
Observe that
\begin{equation*}
  \OpSpaceA{-\infty,d,\E}
  = \bigcap_{N \in \N} \OpSpaceB{-\infty,d,\E}{N}.
\end{equation*}


\begin{lemma} \label{Lempsi-inftysub}
We have
\[ \OpSpaceB{-\infty,d}{N,N} \subset \OpSpaceB{-\infty,d,\E}{N},\]
where $\E$ is the index family on $\T_d$ given by
$\E=(\varnothing, \varnothing, \N_0,\N_0)$.
\end{lemma}
\begin{proof}
Let $A(\lambda) \in  \OpSpaceB{-\infty,d}{N,N}$ and let
$K_{A(\lambda)} = k(\varrho^d \lambda, q)\,\mass'$ be its kernel with all
the properties described in Definition~\ref{SpaceB}. In particular,
the operator $\widetilde A(\lambda)$ with kernel $k(\lambda, q)\,\mass'$
is such that $\lambda^N \widetilde A(\lambda)$ belongs to
$\Sigma^{0}(\overline{\Lambda},\Psi^{-\infty,Nd}_b(M))$.
By definition, we only need to analyze $k(\varrho^d \lambda, q)$ locally in
coordinates near the face $\ffi$.  By symmetry it suffices to consider the
kernel only away from one of the lateral boundaries of $M^2_b$; for
instance, away from the left boundary $\lb$.  Since our kernels are smooth
in $\theta = \lambda/|\lambda|$ and in the variables on the boundary, we
shall omit these variables in what follows.  Thus consider the coordinates
$q=(x,t)\in \U \subset M^2_b$, with $x$ defining $\ff$ and $t = x'/x$
defining $\rb$, see Figure \ref{fig:proj}.  If $\rho = |\lambda|$, then 
for some $\eps > 0$ we can write
\[ k(\lambda, q) = t^{Nd + \eps} g(\rho,x,t), \]
where $g$ is a symbol in $\Sigma^0(\R_+\times \U)$ that can be
expanded in $x$ at $x = 0$ with index set $\N_0$ of order $Nd + \eps$.
In particular, $k$ has a partial asymptotic expansion at $t=0$ with
index set $E_\rb = \varnothing$ of order $Nd + \eps$.

We now lift $k(x^d \lambda, q)$ to $\T_d$. Near $\ff$ and $\ffi$,
the variable $r = \rho^{-1/d}$ defines $\ffi$ and $v = x/r$ defines
$\ff$, and in these coordinates,
\[ k(x^d \lambda, q) = t^{Nd + \eps} g(v^{d},rv,t). \]
The asymptotic properties of $g$ imply that $g(v^{d},rv,t)$ can be
expanded in $r$ and $v$ with index set $\N_0$ of order $Nd +\eps$.
On the other hand, near $\ffi$ and $bi$, $x$ defines $\ffi$ and
$w = v^{-1}$ defines $bi$, and in these coordinates,
\[ k(x^d \lambda, q) = t^{Nd + \eps} g(w^{-d},x,t). \]
Now, since $\lambda^N \widetilde A(\lambda)\in
\Sigma^{0}(\overline{\Lambda},\Psi^{-\infty,Nd}_b(M))$, the function $g$
can actually be written as $g(\rho,x,t)=\rho^{-N} h(\rho,x,t)$, where
$h$ is a symbol in $\Sigma^0(\R_+\times \U)$. Therefore,
\[ k(x^d \lambda, q) = t^{Nd + \eps} w^{Nd} h(w^{-d},x,t). \]
The asymptotic properties of $g$ imply that $h(w^{-d},x,t)$ can be
expanded in $x$ at $x=0$ with index set $\N_0$ of order $Nd + \eps$.
In conclusion, we have proven that $k$ defines a function on $\T_d$
that vanishes to order $Nd$ at $bi$ and has partial expansions of
order $Nd+\eps$ with index sets $E_\rb=\varnothing$, $E_\ff=\N_0$
and $E_\ffi=\N_0$. The same is true with $\rb$ replaced by $\lb$,
thus $A(\lambda)\in \OpSpaceB{-\infty,d,\E}{N}$ with
$\E=(\varnothing,\varnothing,\N_0,\N_0)$.
\end{proof}

\subsection*{Composition of pseudodifferential operators}

In order to prove essential composition properties of our new
parameter-dependent spaces, we need to review how
$b$-pseudodifferential operators are composed. Let $A$ and $B$ be
operators on $C^\infty(M)$ with Schwartz kernels $K_A$ and $K_B$,
respectively, that are smooth on $M^2$ and vanish to infinite order
at the boundary $Y^2=\partial M^2$. Then we know that $AB$ is also a
smoothing operator, and
\begin{equation}\label{KAB}
K_{AB}(u,w) = \int_{v\in M} K_A(u,v)K_B(v,w).
\end{equation}
We can write this purely in terms of pullbacks and pushforwards of
distributions as follows. Let $\pi_F$, $\pi_S$, $\pi_C:M^3 \to M^2$
be the maps
\[ \pi_F(u,v,w) = (u,v),\quad \pi_S(u,v,w)= (v,w),\quad
   \pi_C(u,v,w)= (u,w) \]
($F$, $S$, and $C$ stand for `first', `second', and `composite').
Writing $K_A = k_A\, \mass'$ and $K_B = k_B\, \mass'$, where $k_A$
and $k_B$ are smooth functions on $M^2$ vanishing to infinite order
at the boundary $Y^2$, we have
\[ (\pi_C^*\mass\, \pi_F^* K_A\, \pi_S^*
   K_B )(u,v,w) = k_A(u,v)\,k_B(v,w)\, \mass(u)\mass(v)\mass(w),\]
where on the left-hand side, $\mass$ represents the fixed $b$-measure
on $M$ lifted to $M^2$ under the left projection, that is,
$\mass(u,w) =\mass(u)$ for all $(u,w) \in M^2$. In particular,
$\pi_C^*\mass\, \pi_F^* K_A\, \pi_S^* K_B$ is a density on $M^3$ and so
its pushforward to $M^2$ under $\pi_C$ is well-defined. By \eqref{KAB}
and the definition of the pushforward $(\pi_C)_*$ we get
\begin{equation} \label{KABpf}
\mass\, K_{AB} = (\pi_C)_*(\pi_C^*\mass\, \pi_F^* K_A\, \pi_S^* K_B).
\end{equation}
This identity shows that we can determine the Schwartz kernel of $AB$
by analyzing pullbacks, products, and pushforwards of the kernels $K_A$
and $K_B$. Now,  since our operators are actually in $\Psi^*_b(M)$, 
in order to get a similar identity for the Schwartz kernel of $AB$, 
we first introduce the blown-up manifold $M^3_b$.

\begin{figure} \centering
\setlength{\unitlength}{.25mm} \centering
\begin{picture}(260,110)(-40,-15) { 
\put(5,40){\vector(1,0){40}} \put(160,40){\vector(1,0){40}}}
\put(-110,-10){ \begin{picture}(100,100)(-50,-50) 
\put(0,0){\circle*{4}} 
\put(-90,-3){\small $\begin{array}{ll}\text{Blow up}\\ 
\text{corner first} \end{array}$} 
\put(-25,0){\vector(1,0){18}} \put(25,35){$M^3$}
\put(0,0){\line(-1,-1){20}} \put(0,0){\line(1,0){45}}
\put(0,0){\line(0,1){40}} 
\end{picture} }
\put(30,-10){ \begin{picture}(100,100)(-50,-50) 
\put(-30,-25){\line(-1,-1){20}} \put(20,-7){ \line(1,0){30}}
\put(-9,20){\line(0,1){30}} \qbezier(-9,20)(20,17)(22,-7)
\qbezier(-30,-25)(-36,10)(-9,20) \qbezier(-30,-25)(10,-36)(22,-7)
\put(45,-40){\vector(-1,0){85}}
\put(45,-40){\vector(-2,3){50}} \put(45,-40){\vector(0,1){30}}
\put(42,-55){\small $\begin{array}{ll}\text{Blow up}\\ 
\text{axes next} \end{array}$}
\end{picture} }
\put(200,-10){
\begin{picture}(100,100)(-50,-50) 
\put(-32,-15){\line(-1,-1){20}} \put(-20,-28){\line(-1,-1){20}}
\qbezier(-32,-15)(-22,-15)(-20,-28)
\put(21,0){\line(1,0){30}} \put(20,-15){\line(1,0){33}}
\qbezier(22,0)(12,-7)(20,-15)
\put(3,20){\line(0,1){30}} \put(-15,20){\line(0,1){30}}
\qbezier(3,20)(-2,12)(-15,20) \qbezier(3,20)(18,15)(22,0)
\qbezier(-32,-15)(-36,10)(-15,20)
\qbezier(-20,-28)(10,-36)(20,-15)
\put(-11,-7){\small\emph{ff}} \put(-40,-31){\footnotesize\emph{Sb}}
\put(30,-11){\footnotesize\emph{Cb}} \put(-13,32){\footnotesize\emph{Fb}}
\put(0,-42){\small\emph{rb}} \put(22,17){\small\emph{lb}}
\put(-50,10){\small\emph{mb}} \put(30,40){$M^3_b$}
\end{picture} } 
\end{picture}
\caption{The blown-up manifold $M^3_b$ and its boundary
hypersurfaces.} \label{fig:M3b}
\end{figure}
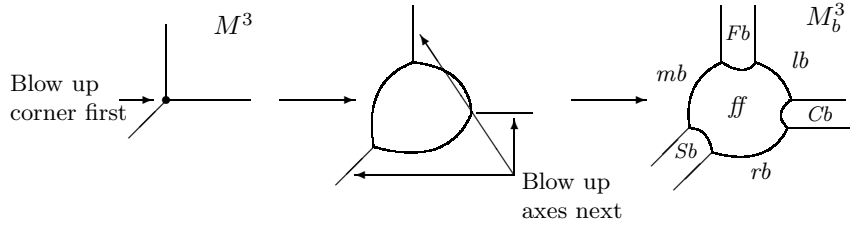

The manifold $M^3_b$ is defined by blowing-up (that is, introducing polar
coordinates around) the manifold $Y^3$ in $M^3$ and then blowing-up the
submanifolds coming from the codimension two corners of $M^3$. The
manifold $M^3_b$ along with its various faces are shown in Figure
\ref{fig:M3b}.
Let $\pi_{F,b}$, $\pi_{S,b}$, $\pi_{C,b}:M^3_b \rightarrow M^2_b$ be
the maps $\pi_F$, $\pi_S$, $\pi_C$ expressed in the
polar coordinates of $M^3_b$ and $M^2_b$. Then we can express the
composition \eqref{KABpf} in terms of these new maps:
\begin{equation} \label{KABbpf}
\mass\, K_{AB} = (\pi_{C,b})_*(\pi_{C,b}^*\mass\, \pi_{F,b}^*
K_A\, \pi_{S,b}^* K_B).
\end{equation}
Written in this way, $\mass$, $K_A$, $K_B$, and $K_{AB}$ are understood
to be lifted to $M^2_b$. The formula \eqref{KABbpf} is the key to proving
composition properties of our parameter-dependent operators.

\subsection*{Composition theorems for parameter-dependent operators}

We begin by stating a composition result whose proof is almost exactly the 
same as the proof of \cite[Theorem~4.4]{LoRes01}, so we leave out the
details. 

\begin{theorem} \label{thm:comp1} We have
\[ x^\nu \OpSpaceA{\mu,p,d}\circ x^{\nu'}\OpSpaceA{\mu'\!,p'\!,d}
  \subset x^{\nu+\nu'} \OpSpaceA{\mu+\mu'\!,p+p'\!,d}. \]
\end{theorem}

The following theorem is proved in \cite[Proposition~4.1]{LoRes201}.

\begin{theorem} \label{thm:comp2}
We have
\begin{gather*}
x^\nu \Psi^{\mu}_{b}(M)\circ x^{\nu'}\OpSpaceA{\mu'\!,p,d}\subset
 x^{\nu+\nu'} \OpSpaceA{\mu+\mu'\!,p,d};
\intertext{and}
x^{\nu'}\OpSpaceA{\mu'\!,p,d}\circ x^\nu \Psi^{\mu}_{b}(M) \subset
 x^{\nu+\nu'} \OpSpaceA{\mu+\mu'\!,p,d}.
\end{gather*}
\end{theorem}

This theorem states that the parameter-dependent spaces
$\OpSpaceA{*,p,d}$ are closed under composition with
cone pseudodifferential operators. We next consider how these
spaces behave under composition with the calculus with bounds and
their parameter-dependent versions. The next theorem is
established by following the proof of \cite[Proposition~4.1]{LoRes201},
taking into account the bounds. To avoid reproducing the proof of
loc.~cit., we shall omit the details.

\begin{theorem} \label{thm:comp3}
The spaces
$\OpSpaceB{-\infty,d}{m,N}$ and $\OpSpaceB{-\infty,d,\E}{N}$ for
any $m, N \in \N$ are closed under compositions with
$\Psi^m_b(M)$.  Let $p/d \in - \N$. Then for any $N > 0$, we have
\begin{gather*}
\Psi^{-\infty,\,2Nd}_{b}(M)\circ \OpSpaceA{\mu,p,d}
\subset \OpSpaceB{-\infty,d}{m,N},
\\
\OpSpaceA{\mu,p,d}\circ \Psi^{-\infty,\,2Nd}_{b}(M)
\subset \OpSpaceB{-\infty,d}{m,N},
\intertext{where $m = - p/d$;}
\OpSpaceB{-\infty,d}{m',N}\circ \OpSpaceA{\mu,p,d}
\subset \OpSpaceB{-\infty,d}{m,N},
\\
\OpSpaceA{\mu,p,d}\circ \OpSpaceB{-\infty,d}{m',N}
\subset \OpSpaceB{-\infty,d}{m,N},
\end{gather*}
where $m = \min\{m',-p/d\}$. Finally, the space
$\OpSpaceB{-\infty,d,\E}{N}$ is closed under composition with
$\OpSpaceA{\mu,p,d}$, for instance,
\[ \OpSpaceA{\mu,p,d}\circ \OpSpaceB{-\infty,d,\E}{N}
\subset \OpSpaceB{-\infty,d,\E}{N}. \]
\end{theorem}

We next consider composition in our first new parameter-dependent calculus.

\begin{theorem} \label{thm-compositemN}
We have
\[ \OpSpaceB{-\infty,d}{m,N_1} \circ
 \OpSpaceB{-\infty,d}{m',N_2} \subset \OpSpaceB{-\infty,d}{m + m', N}, \]
where $N = \min\{N_1,N_2\}$.
\end{theorem}
\begin{proof}
Since $\OpSpaceB{-\infty,d}{m, N_1} \subset \OpSpaceB{-\infty,d}{m, N}$
and $\OpSpaceB{-\infty,d}{m', N_2} \subset \OpSpaceB{-\infty,d}{m', N}$,
which follows from the definition of these spaces, we
may assume that $N = N_1=N_2$. Thus given $A \in
\OpSpaceB{-\infty,d}{m,N}$ and $B \in \OpSpaceB{-\infty,d}{m',N}$,
we need to show that $A B \in \OpSpaceB{-\infty,d}{m + m',N}$. To
simplify the notation, we assume that $M = [0,1)_x$ and $\Lambda =\R_+$.
The argument in the general case is basically the same, the main difference 
being the appearance of the tangential variables on $Y$ that make the proof 
notationally more complicated. We will use projective coordinates (see 
Figure \ref{fig:proj}).

In the following, we denote by $x$, $x'$, $x''$ the coordinates on
the left, middle, and right factors of $M^3 = [0,1)^3$ and we assume
that $\mass = |dx/x|$. To show that $A B \in \OpSpaceB{-\infty,d}{m +
m',N}$, we use the formula \eqref{KABbpf} above.
To do so, we will assume that the lifted kernel
$\pi_{C,b}^* \mass \,\pi_{F,b}^* K_A \, \pi_{S,b}^*K_B$
is supported in a neighborhood $\U\subset M^3_b$ of the intersection
of $mb$, $\ff$, and $Fb$ (see Figure \ref{fig:M3b}).
On $\U$ we introduce projective coordinates as follows.
First, we blow-up the origin in $M^3$ and define,
away from the hypersurface $\{x''=0\}$, the coordinates
$(s,s',x'')$ with $s=x/x''$ and $s'= x'/x''$.
Next, we blow-up the $x''$-axis to get $M_b^3$, and define on $\U$
the coordinates $(s,t,x'')$ with  $t=s'/s = x'/x$. In conclusion,
we get the projective coordinates
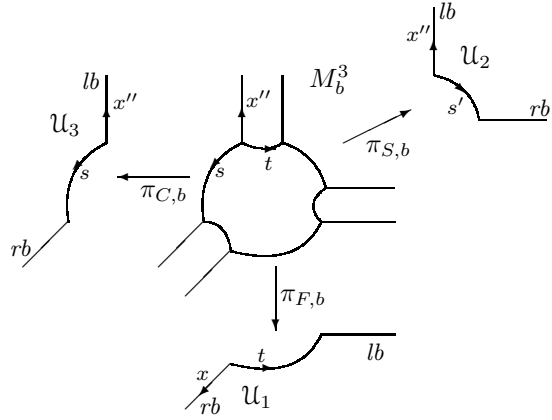
\begin{figure} \centering
\setlength{\unitlength}{.30mm} \centering 
\begin{picture}(260,180)(-130,-100)
\put(-32,-15){\line(-1,-1){20}} \put(-20,-28){\line(-1,-1){20}}
\qbezier(-32,-15)(-22,-15)(-20,-28)
\put(19,0){ \line(1,0){30}} \put(20,-15){\line(1,0){33}}
\qbezier(22,0)(12,-7)(20,-15)
\put(3,20){\line(0,1){30}} \put(-15,20){\line(0,1){30}}
\qbezier(3,20)(-5,15)(-15,20) \qbezier(3,20)(18,15)(22,0)
\qbezier(-32,-15)(-36,10)(-15,20) \qbezier(-20,-28)(10,-36)(20,-15)
\put(15,45){$M^3_b$}
\put(-15,35){\vector(0,1){1}} 
\put(-1,17.5){\vector(1,0){1}} 
\put(-28,10){\vector(-1,-1){1}}
\put(-13,36){\footnotesize $x''$} 
\put(-27,5){\footnotesize $s$}
\put(-5,8){\footnotesize $t$} 
%
\put(-20,-78){\line(-1,-1){20}} \put(20,-65){\line(1,0){33}}
\qbezier(-20,-78)(10,-86)(20,-65) 
\put(-34,-100){\small \emph{rb}} \put(40,-76){\small \emph{lb}}
\put(-15,-96){$\U_1$} \put(0,-35){\vector(0,-1){28}}
\put(2,-50){$\pi_{F,b}$}
\put(-2,-80){\vector(1,0){1}} 
\put(-33,-91){\vector(-1,-1){1}}
\put(-35,-85){\footnotesize $x$} 
\put(-8,-77){\footnotesize $t$}
\put(70,50){\line(0,1){30}} \put(90,30){\line(1,0){30}}
\qbezier(70,50)(87,46)(90,30) 
\put(112,32){\small \emph{rb}} \put(72,75){\small \emph{lb}}
\put(82,55){$\U_2$} \put(30,20){\vector(2,1){30}}
\put(40,16){$\pi_{S,b}$}
\put(70,65){\vector(0,1){1}} 
\put(83,43){\vector(1,-1){1}} 
\put(58,65){\footnotesize $x''$} 
\put(76,33){\footnotesize $s'$} 
\put(-92,-15){\line(-1,-1){20}} \put(-75,20){\line(0,1){30}}
\qbezier(-92,-15)(-96,10)(-75,20) 
\put(-118,-30){\small\emph{rb}} \put(-86,44){\small \emph{lb}}
\put(-100,25){$\U_3$} \put(-38,5){\vector(-1,0){32}}
\put(-60,-4){$\pi_{C,b}$} 
\put(-75,35){\vector(0,1){1}} 
\put(-89,9){\vector(-1,-1){1}}
\put(-72,35){\footnotesize $x''$} 
\put(-87,4){\footnotesize $s$}
\end{picture}
\caption{Projective coordinates and projections on $\U$.} \label{fig:Ucoord}
\end{figure}
\begin{equation}\label{Ucoord}
 (s,t,x'')\in \U \quad \text{with }\; s=\frac{x}{x''}
 \;\text{ and }\; t=\frac{x'}{x}.
\end{equation}
By definition, $\pi_{F,b}(s,t,x'')$ is the image of
$\pi_{F}(x,x',x'')=(x,x')$ written in coordinates $(x,x'/x)$ on
$M_b^2$. Similarly, $\pi_{S,b}(s,t,x'')$ is the image of
$\pi_{S}(x,x',x'')=(x',x'')$ written in coordinates $(x'/x'',x'')$.
The appropriate choice of projective coordinates on $M^2_b$ for the
images of $\pi_F$ and $\pi_S$ is illustrated in Figure~\ref{fig:Ucoord}.
By means of \eqref{Ucoord} we finally get
\begin{equation} \label{mbfffs2}
\pi_{F,b}(s,t,x'')=(sx'',t), \quad
\pi_{S,b}(s,t,x'')=(st,x'').
\end{equation}

Let $\U_1 = \pi_{F,b}(\U)$. In the coordinates $(x,t)\in \U_1$
(see Figure~\ref{fig:Ucoord}), the kernel of $A$ is of the form
$K_A = k_1(x^d\lambda,x,t)\,|dx'/x'|$, where for some $\eps > 0$,
$t^{-Nd - \eps} k_1(\lambda,x,t)$ is a symbol in
$\Sigma^0(\R_+\times \U_1)$ that can be expanded at $x=0$ with index
set $\N_0$ of order $Nd + \eps$. Moreover, for each $N'\leq N$,
\begin{equation} \label{expAlambdaxt}
 k_1(\lambda,x,t) = \sum_{j = m}^{N'-1} \lambda^{-j}\, f_{j}(x,t)
 + \lambda^{-N'} k_{1,N'}(\lambda,x,t),
\end{equation}
where all the coefficients satisfy the properties listed in
Definition~\ref{SpaceB}.

Let $\U_2 = \pi_{S,b}(\U)$. In the coordinates $(s',x'')\in
\U_2$ (Figure~\ref{fig:Ucoord}), the kernel of $B$ is of the
form $K_B = k_2((x'')^d\lambda,s',x'')\,|dx''/x''|$, where for some
$\eps > 0$, $(s')^{-Nd-\eps}k_2(\lambda,s',x'')$  is a symbol in
$\Sigma^0(\R_+\times \U_2)$ that can be expanded at $x''=0$
with index set $\N_0$ of order $Nd + \eps$. The function
$k_2(\lambda,s',x'')$ also admits an expansion similar to
\eqref{expAlambdaxt} with the obvious change of variables.
Using the formulas for $\pi_{S,b}$ and $\pi_{F,b}$ in \eqref{mbfffs2},
it follows that on $\U$,
\[ \pi_{C,b}^* \mass \,\pi_{F,b}^* K_A \, \pi_{S,b}^*K_B
    = k_1((sx'')^d\lambda,sx'',t) \, k_2((x'')^d\lambda,st,x'')
      \Big|\frac{dsdtdx''}{stx''}\Big|. \]
Furthermore,
$\pi_{C,b}(s,t,x'')$ is the image of $\pi_C(x,x',x'')=(x,x'')$ written in
coordinates $(x/x'',x'')$ on $M^2_b$, thus
\[ \pi_{C,b}(s,t,x'') = (s,x''). \]
By the definition of pushforward,
\begin{align*}
 (\pi_{C,b})_*(\pi_{C,b}^* \mass \,\pi_{F,b}^* K_A \, \pi_{S,b}^*K_B)
 = k_3((x'')^d\lambda,s,x'') \Big|\frac{dsdx''}{sx''}\Big|,
\end{align*}
where
\begin{equation*}
 k_3(\lambda,s,x'') =
 \int k_1(s^d\lambda,sx'',t)\, k_2(\lambda,st,x'')\frac{dt}{t}
\end{equation*}
From the properties of $A$ and $B$ it follows easily that
$s^{-Nd - \eps} k_3(\lambda,s,x'')$ is a symbol in $\Sigma^0(\R_+\times
\U_3)$ having a partial expansion at $x'' = 0$ with index set $\N_0$
of order $Nd + \eps$. Moreover, the formula \eqref{expAlambdaxt}
corresponding to  $k_2$ (denoting the coefficients by $g_j$ instead
of $f_j$) implies that given $N'\le N$,
\begin{multline*}
k_3(\lambda,s,x'') = \sum_{j=m'}^{N' - 1} \lambda^{-j}\,
\int k_1(s^d\lambda,sx'',t)\, g_j(st,x'')\, \frac{dt}{t} \\
+ \lambda^{-N'} \int k_1(s^d\lambda,sx'',t)\, k_{2,N'}(\lambda,st,x'')\,
\frac{dt}{t}.
\end{multline*}
Now for each $j$, expanding $k_1(\lambda,x,t)$ in $\lambda$ up to order
$N' -  j$, we find that
\begin{equation}\label{expClambda}
 k_3(\lambda,s,x'') = \sum_{j=m+m'}^{N' - 1} \lambda^{-j}\, h_j(s,x'')
 + \lambda^{-N'} k_{3,N'}(\lambda,s,x''),
\end{equation}
where
\begin{equation*}
  h_j(s,x'') =  \sum_{\ell=m'}^{j-m}
  s^{-(j-\ell)d}\int f_{j-\ell}(sx'',t)\, g_{\ell}(st,x'') \frac{dt}{t}
\end{equation*}
and
\begin{multline*}
k_{3,N'}(\lambda,s,x'') = \sum_{j=m'}^{N'-1} s^{-(N'-j)d}
\int k_{1,N'-j}(s^d \lambda,s,x'')\, g_{\ell}(st,x'') \frac{dt}{t} \\
+ \int k_1(s^d\lambda,sx'',t)\, k_{2,N'}(\lambda,st,x'') \frac{dt}{t}.
\end{multline*}
It remains to verify that the coefficients in the expansion
\eqref{expClambda} have the properties required in part (b) of
Definition~\ref{SpaceB}. For the $h_j$'s this follows
from the corresponding properties of the functions $f_{j-\ell}$ and
$g_\ell$. In particular, expanding these functions at $x''=0$ according
to \eqref{expprop} we get the required expansion for $h_j$.
Notice that these expansions are partial expansions with index set $\N_0$
of order $2Nd-jd+\ell d$ for $f_{j-\ell}$ and $2Nd-\ell d$ for $g_\ell$,
which are both of order greater than $2Nd - jd$. Therefore, the resulting
asymptotic expansion for $h_j(s,x'')$ at $x''=0$ is of the same type.
The properties of $k_{3,N'}(\lambda,s,x'')$ follow in a similar manner.
\end{proof}

\begin{remark}\label{neighborhoods}
In the previous proof, we restricted ourselves to a neighborhood
$\U\subset M^3_b=[0,1)^3_b$ of the intersection of the faces $mb$, $\ff$,
and $Fb$ (see Figure~\ref{fig:M3b}). But in fact, in this model case,
we need six coordinate patches to cover the entire manifold.
However, by symmetry, there are only three patches that require slightly 
different treatments. For instance, we could choose in addition to $\U$
a neighborhood $\mathcal V$ of the intersection of $Fb$, $\ff$, and $\lb$,
and a neighborhood $\mathcal W$ of the intersection of $\ff$, $Cb$,
and $\lb$, to complete a set of representative local coordinates,
see Figure~\ref{fig:coordM3b}.
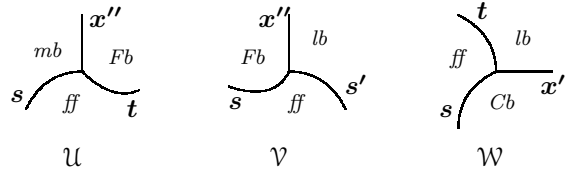
\begin{figure} \centering
\centering \setlength{\unitlength}{.5mm}
\begin{picture}(120,45)(-50,-10) 
\put(-45,10){\begin{picture}(100,80) 
\qbezier(0,0)(-10,0)(-15,-10)
\qbezier(0,0)(8,-8)(15,-5)\put(0,0){\line(0,1){15}}
\put(-19,-8){\boldmath $s$}\put(12,-12){\boldmath $t$} 
\put(2,12){\boldmath $x''$}
\put(-13,4){\footnotesize \emph{mb}} \put(-5,-10){\footnotesize \emph{ff}}
\put(7,2){\footnotesize \emph{Fb}} \end{picture} }
\put(-50,-15){$\mathcal U$}
\put(5,-15){$\mathcal V$}
\put(60,-15){$\mathcal W$}

\put(10,10){\begin{picture}(100,80) 
\put(0,0){\line(0,1){15}} \qbezier(0,0)(10,0)(15,-10)
\qbezier(0,0)(-4,-8)(-16,-4) \put(-16,-10){\boldmath $s$}
\put(15,-6){\boldmath $s'$}
\put(-8,12){\boldmath $x''$} \put(-13,2){\footnotesize \emph{Fb}}
\put(0,-10){\footnotesize \emph{ff}}
\put(6,7){\footnotesize \emph{lb}} \end{picture} }

\put(65,10){\begin{picture}(100,80) 
\qbezier(0,0)(0,10)(-10,15) \qbezier(0,0)(-10,-5)(-10,-15)
\put(0,0){\line(1,0){15}} \put(-15,-12){\boldmath $s$}
\put(12,-6){\boldmath $x'$} \put(-5,14){\boldmath $t$} 
\put(-12,2){\footnotesize \emph{ff}}
\put(-2,-10){\footnotesize \emph{Cb}}
\put(5,7){\footnotesize \emph{lb}} \end{picture} }

\end{picture}
\caption{Three representative coordinate patches on $M^3_b$.}
\label{fig:coordM3b}
\end{figure}
Since the calculations on $\mathcal V$ and $\mathcal W$ are similar in
nature, and in order to avoid an overloading of technical computations,
we decided to skip them. Nonetheless, to demonstrate these different
treatments without repeating our arguments, in the proofs of
Theorem~\ref{thm-full} and Theorem~\ref{thm-boundsasym}
we will work on $\mathcal V$ and $\mathcal W$, respectively.
\end{remark}

For index sets $E$ and $F$, we define the extended union of these sets by 
\begin{equation} 
E \overline{\cup} F = E \cup F \cup \{(z,k+\ell + 1) \st 
(z,k) \in E,\ (z,\ell) \in F\}.
\label{extunion}
\end{equation}
Given $\E=(E_{lb},E_{rb},E_\ff,E_\ffi)$ and 
$\F= (F_{lb}, F_{rb},F_\ff,F_\ffi)$, we define the index family $\E
\hat{\circ} \F = (G_{lb},G_{rb},G_\ff,G_\ffi)$ as follows:
\begin{equation}
\begin{array}{c} G_{lb}= E_{lb}\overline{\cup}(E_\ff + F_{lb}),\;
G_{rb}=(E_{rb} + F_\ff) \overline{\cup}F_{rb}, \\ 
G_\ff=(E_\ff + F_\ff)\overline{\cup}(E_{lb}+ F_{rb}),\; \text{ and }\;
G_\ffi= E_\ffi+F_\ffi. 
\end{array} \label{compindexset}
\end{equation}

Our second new parameter-dependent space has the following properties.
\begin{theorem}  \label{thm-full}
Provided that $E_{rb} + F_{lb} > 0$, we have
\[ \OpSpaceB{-\infty,d,\E}{N}\circ \OpSpaceB{-\infty,d,\F}{N}
\subset \OpSpaceB{-\infty,d, \E \hat{\circ}\F }{N}.\]
\end{theorem}
\begin{proof} Let $A \in \OpSpaceB{-\infty,d,\E}{N}$ and
$B \in \OpSpaceB{-\infty,d,\F}{N}$. We will use the formula
\eqref{KABbpf} to show that $A B \in 
\OpSpaceB{-\infty,d, \E\hat{\circ} \F }{N}$. As in the previous proof, we 
assume that $M = [0,1)_x$ and $\Lambda = \R_+$. We also use projective 
coordinates near $\lb$ on the product $[0,1)^2_b$ (see Figure \ref{fig:proj}).

Let $x$, $x'$, $x''$ be the coordinates on the left, middle, and
right factors of $[0,1)^3$ and assume that  $\mass = |dx/x|$. 
In this proof we now assume the lifted kernel 
$\pi_{C,b}^* \mass \,\pi_{F,b}^* K_A \,\pi_{S,b}^*K_B$ to be supported 
in a neighborhood $\mathcal V\subset M^3_b$ of the intersection of $Fb$,
$\ff$, and $\lb$ (see Figure \ref{fig:M3b}). On $\mathcal V$ we may use the 
coordinates 
\begin{equation*}
(s,s',x'')\quad\text{with } s=\frac{x}{x'}\;\text{ and }\; s'=\frac{x'}{x''}
\quad \text{(see Figure \ref{fig:coordM3b})}.
\label{coordsnearlb}
\end{equation*}
The projections $\pi_{F,b}$, $\pi_{S,b}$, and $\pi_{C,b}$ all map
$\mathcal V$ onto a neighborhood of $\lb$ in $M^2_b$, and we have
\begin{equation} \label{fsfflb}
\begin{split} 
 \pi_{F,b}(s,s',x'') &= (s,s'x''),\quad\pi_{S,b}(s,s',x'')=(s',x''),\\ 
 & \pi_{C,b}(s,s',x'')=(ss',x'').
\end{split}
\end{equation}
If $r = \lambda^{-1/d}$ and $v'=x'/r$, then near $\lb$ in $\T_d$, we can 
write $K_A=k_1(r,s,v')|dx'/x'|$, where for some $\eps> 0$, $k_1(r,s,v')$ 
has expansions at $r=0$, $v'=0$, and $s=0$, with index sets $E_\ffi$, 
$E_\ff$, and $E_\lb$ of order $Nd + \eps$, respectively, and for $v' \geq 1$, 
$k_1(r,s,v')= (v')^{-Nd - \eps} \widetilde{k_1}(r,s,v')$ where
$\widetilde{k_1}(r,s,v')$ is a symbol of order $0$ in $v'$ and has
expansions at $r=0$ and $s=0$ with index sets $E_\ffi$ and $E_\lb$ of order 
$Nd + \eps$, respectively. Similarly, $K_B=k_2(r,s,v') |dx'/x'|$ where 
$k_2(r,s,v')$ has analogous properties as $k_1(r,s,v')$ but with index sets 
given by $\F$. Using the formulas in \eqref{fsfflb}, it follows that
\[ \pi_{C,b}^* \mass \,\pi_{F,b}^* K_A \,
\pi_{S,b}^* K_B = k_1(r,s,s'x''/r) k_2(r,s',x''/r)
\Big|\frac{dsds'dx''}{ss'x''}\Big|.\] 
Hence, as $\pi_{C,b}(s,s',x'')=(ss',x'')$, by the definition of 
pushforward, we obtain
\[ (\pi_{C,b})_*(\pi_{C,b}^* \mass \,\pi_{F,b}^* K_A \, \pi_{S,b}^* K_B) 
   = k_3(r,s,v')\Big|\frac{dsdx'}{sx'}\Big|, \]  
where 
\[ k_3(r,s,v')=\int k_1(r,s/s',s'v') k_2(r,s',v')\frac{ds'}{s'}. \] 
Now the asymptotic properties of $k_1$ and $k_2$ together with Melrose's 
pushforward theorem (see \cite[Appendix]{LoRes01}) imply that
$k_3(r,s,v')$ has expansions at $r=0$, $s=0$, and $v'=0$, with index sets 
$E_\ffi+F_\ffi$, $E_\ff+F_\ff$, and $E_\lb\overline{\cup}(E_\ff + F_\lb)$ 
of order $Nd + \eps$, respectively. Moreover, for $v' \geq 1$, 
$k_3(r,s,v')= (v')^{-Nd - \eps} \widetilde{k_3}(r,s,v')$ where 
$\widetilde{k_3}(r,s,v')$ is a symbol of order $0$ in $v'$ and has
expansions at $r=0$ and $s=0$ with index sets $E_\ffi+F_\ffi$ and
$E_\lb\overline{\cup}(E_\ff + F_\lb)$ of order $Nd+\eps$, respectively.
\end{proof}

Finally, we consider the composition of our two new parameter-dependent 
spaces.

\begin{theorem}  \label{thm-boundsasym}
We have
\begin{gather*}
\OpSpaceB{-\infty,d}{m,N} \circ \OpSpaceB{-\infty,d,\E}{N}
\subset \OpSpaceB{-\infty,d, \E}{N},\\
\OpSpaceB{-\infty,d,\E}{N} \circ \OpSpaceB{-\infty,d}{m,N} 
\subset \OpSpaceB{-\infty,d, \E}{N}.
\end{gather*}
\end{theorem}
\begin{proof}
Let $A \in \OpSpaceB{-\infty,d}{m,N}$ and $B\in \OpSpaceB{-\infty,d,\E}{N}$. 
As in our previous proofs, we will use the composition formula \eqref{KABbpf} 
to show that $A B \in \OpSpaceB{-\infty,d, \E}{N}$. Again, we consider 
$M = [0,1)_x$, $\Lambda = \R_+$, and introduce the following coordinates 
on $[0,1)^2_b$  (see Figure \ref{fig:proj}):
\begin{equation} \label{coordslb3}
(s,x') \text{ near $\lb$, and } (x,t) \text{ near $\rb$, where } 
s=\frac{x}{x'} \text{ and } t=\frac{x'}{x}.
\end{equation}

Let $x$, $x'$, $x''$ be the coordinates on the left, middle, and
right factors of $[0,1)^3$ and assume that  $\mass = |dx/x|$. 
According to Remark~\ref{neighborhoods}, in the present proof we will assume 
that the lifted kernel $\pi_{C,b}^*\mass\,\pi_{F,b}^* K_A\,\pi_{S,b}^*K_B$ 
is supported in a neighborhood $\mathcal W$ of the intersection of
$\ff$, $Cb$, and $\lb$ in $M^3_b$ (see Figure \ref{fig:M3b}).
Here, we may use the coordinates $(s,x',t)\in\mathcal W$, where $s=x/x''$ 
and $t=x''/x'$ (see Figure \ref{fig:coordM3b}). The projections $\pi_{C,b}$ 
and $\pi_{F,b}$ map $\mathcal W$ onto a neighborhood of $\lb$ in
$M^2_b$, and $\pi_{S,b}$ maps $\mathcal W$ onto a neighborhood of
$\rb$. Moreover, in the coordinates \eqref{coordslb3} on $M^2_b$
near $\lb$, we have
\begin{equation} 
\pi_{F,b}(s,x',t)=(st,x'),\quad \pi_{C,b}(s,x',t)=(s,x't),
\label{ffcslb3}
\end{equation}
and near $\rb$, we have
\begin{equation} \label{ffcslb4}
  \pi_{S,b}(s,x',t)=(x',t).
\end{equation} 
Near $\lb$ in $M^2_b$, we can write $K_A = k_1((x')^d\lambda,s,x')|dx'/x'|$, 
where for some $\eps > 0$, $s^{-Nd - \eps} k_1(\lambda,s,x')$ is a
symbol of order $0$ in all variables that can be expanded at $x' = 0$ with 
index set $\N_0$ of order $Nd + \eps$, and for each $N' \leq N$, $k_1$ can 
be written in the form 
\begin{equation}\label{expAlambdasx'} 
k_1(\lambda,s,x') = \sum_{j = m}^{N'-1}
\lambda^{-j}\, f_j(s,x') + \lambda^{-N'} f_{N'}(\lambda,s,x'),
\end{equation} 
where $f_{N'}(\lambda,s,x')$ satisfies the same conditions as 
$k_1(\lambda,s,x')$ but with $N$ replaced by $2Nd -N'd$, and for each 
$m \leq j \leq N' - 1$, $f_j(s,x')$ satisfies the same conditions as 
$k_1(\lambda,s,x')$ but with $Nd$ replaced by $2Nd - j d$; of course, 
without the condition \eqref{expAlambdasx'}. If $v=x/r$ where 
$r = \lambda^{-1/d}$, then near $\rb$ in $\T_{d}$, we can write 
$K_B=k_2(r,v,t) |dx'/x'|$, where for some $\eps > 0$, $k_2(r,v,t)$ has 
expansions at $r=0$, $v =0$, and $t=0$, with index sets $E_\ffi$,
$E_\ff$, and $E_\rb$ of order $Nd + \eps$, respectively, and for $v \geq 1$, 
$k_2(r,v,t) = v^{-Nd -\eps} \widetilde{k_2}(r,v,t)$ where 
$\widetilde{k_2}(r,v,t)$ is a symbol of order $0$ in $v$ and has expansions 
at $r=0$ and $t=0$ with index sets $E_\ffi$ and $E_\rb$ of order $Nd +\eps$, 
respectively. Using the formulas for $\pi_{F,b}$ and $\pi_{S,b}$ in 
\eqref{ffcslb3} and \eqref{ffcslb4}, it follows that on $\mathcal W
\subset M^3_b$, 
\[ \begin{split} 
\pi_{C,b}^* \mass \,\pi_{F,b}^* K_A \, \pi_{S,b}^* K_B 
& = k_1((x')^d\lambda,st,x') k_2(r,x'/r,t)\Big|\frac{dsdtdx'}{stx'} \Big| \\ 
& = k_1((x'/r)^d,st,x') k_2(r,x'/r,t) \Big|\frac{dsdtdx'}{stx'} \Big|, 
\end{split} \] 
since $r=\lambda^{-1/d}$. Hence, as $\pi_{C,b}(s,x',t)=(s,x't)$, working
out the definition of pushforward we obtain 
\[ (\pi_{C,b})_*(\pi_{C,b}^* \mass \, \pi_{F,b}^* K_A \,\pi_{S,b}^* K_B) 
 = k_3(r,v',s) \Big| \frac{ds dx'}{sx'}\Big|, \] 
where 
\[ k_3(r,v',s)=\int k_1((v'/t)^d,st,rv') k_2(r,v'/t,t) \frac{dt}{t}. \]
Using the asymptotic properties of $k_1(\lambda,s,x')$ given in
\eqref{expAlambdasx'} and the asymptotic properties of $k_2(r,v,t)$,
one can show that $k_3(r,v',s)$ has expansions at $r=0$, $v'=0$, and
$s=0$, with index sets $F_\ffi$, $F_\ff$, and $\varnothing$ of order 
$Nd + \eps$, respectively, and for $v' \geq 1$, $k_3(r,v',s) = 
(v')^{-Nd - \eps} \widetilde{k_3}(r,v',s)$ where $\widetilde{k_3}(r,v',s)$ 
is a symbol of order $0$ in $v'$ and has expansions at $r=0$ and $s=0$ with 
index sets $F_\ffi$ and $\varnothing$ of order $Nd + \eps$, respectively.
\end{proof}

\section{Resolvents and parametrix construction}
\label{Parametrix}

In this section we let $\mu>0$ and consider a cone pseudodifferential 
operator $A\in x^{-\mu}\Psi_b^{\mu}(M)$, where $x$ is a boundary 
defining function for $\partial M$ and $\Psi_b^{\mu}(M)$ is the class 
of $b$-operators, cf. Definition~\ref{b-Operators}. 
It is well-known (see e.g. \cite{RBM2} or \cite{Sz91}) that $A$ can be 
extended as a bounded operator 
\begin{equation}\label{BoundedOp}
  A:x^{\alpha}H_b^s(M)\to x^{\alpha-\mu}H_b^{s-\mu}(M),
\end{equation}
where the space $H_b^s(M)$ is defined as follows.  We fix a $b$-measure
$\mass$ and let $L_b^2(M)$ be the Hilbert space of square integrable 
functions with respect to $\mass$. For $s\in\N$, the space $H_b^s(M)$
consists of all $u\in L_b^2(M)$ such that $Pu\in L_b^2(M)$ for every 
differential operator $P\in \Diff_b^s(M)$. For an arbitrary $s\in\R$, 
the space $H_b^s(M)$ can be defined by duality and interpolation. 

\begin{remark}\label{embeddings}
For $s\ge s'$ and $\alpha\ge \alpha'$ the embedding
$x^{\alpha}H_b^s(M)\embed x^{\alpha'}H_b^{s'}(M)$ is continuous. If
$\alpha> \alpha'$, then it is compact if $s>s'$, Hilbert-Schmidt if
$s>s'+\frac{n}{2}$, and trace class if $s>s'+n$, where $n=\dim M$.
\end{remark}

For $A\in x^{-\mu}\Psi_b^{\mu}(M)$ we let $\sym_\mu(A)$ be the totally
characteristic principal symbol of $x^\mu A$ in $\Psi_b^{\mu}(M)$. The 
operator $A$ is said to be $b$-elliptic if $\sym_\mu(A)$ is invertible on 
${}^bT^*M\setminus 0$, where ${}^bT^*M$ denotes the $b$-cotangent bundle, 
cf. \cite{RBM2}. The Fredholm property of \eqref{BoundedOp} is determined 
by the indicial family (or conormal symbol) $\widehat{A}(z)$ associated 
with $A$. It is defined as the operator family 
\begin{equation*}
  \widehat{A}(z): C^\infty(Y)\to C^\infty(Y): 
  u\mapsto x^{\mu-z}A(x^z \tilde u)|_{x=0},
\end{equation*}
where $Y=\partial M$ and $\tilde u$ is some extension of $u$. The set
\begin{equation*}
  \spec_b(A)=\{z\in\C\st \widehat{A}(z):H^\mu(Y)\to L^2(Y) 
  \text{ is not invertible}\}	
\end{equation*}
is called the \emph{boundary spectrum of} $A$. If $A$ is $b$-elliptic, 
then its boundary spectrum is discrete and we have the following result, 
cf. \cite{RBM2, MM83, Sz91}.

\begin{theorem}
If $A\in x^{-\mu}\Psi_b^{\mu}(M)$ is $b$-elliptic, then for every
$\alpha\in\R$ such that  $\mathrm{spec}_b(A)\cap \{z\in\C \st 
\Im z=-\alpha\}=\varnothing$, the operator \eqref{BoundedOp}
is Fredholm for every $s\in\R$.
\end{theorem}

In order to ensure the existence of the resolvent and be able to describe 
it within our calculus, we need a natural notion of parameter-dependent 
ellipticity that resembles Agmon's condition at the symbol level and takes 
into account the singular global behavior of the operator near the boundary. 
Following \cite{GiHeat01} we will define the parameter-ellipticity with 
help of a model operator $A_\wedge$ living on the model cone
$Y^\wedge:=\overline{\R}_+\times Y$. More precisely, with 
$A$ we associate the operator 
\begin{equation}\label{ModelOperator}
 A_\wedge:C_c^\infty(Y^\wedge)\to C^\infty(Y^\wedge): 
 u\mapsto  \lim_{\varrho\to 0}\varrho^\mu\kappa_\varrho 
 \varphi A(\psi\kappa_\varrho^{-1} u),
\end{equation}
where $\kappa_\varrho$ is defined by 
$(\kappa_\varrho u)(x,y):= u(\varrho x,y)$, $\varrho>0$, and where $\varphi$ 
and $\psi$ are smooth functions supported in a collar neighborhood of $Y$ 
($=\partial M=\partial Y^\wedge$) so that $\psi\kappa_\varrho^{-1} u$ and
$\varphi A(\psi\kappa_\varrho^{-1} u)$ can be regarded as functions on 
both manifolds $M$ and $Y^\wedge$. 

On $Y^\wedge$ it is convenient to introduce Schulze's (cone) Sobolev spaces
$\mathcal{K}^{s,\alpha}(Y^\wedge)$ for $s,\alpha\in\R$, defined as follows. 
Let $\omega\in C^\infty_c(\overline\R_+)$ with $\omega(r)=1$ near $r=0$. 
Then the space $\mathcal{K}^{s,\alpha}(Y^\wedge)$ consists of distributions 
$u$ such that $\omega u\in r^\alpha H^s_b(Y^\wedge)$, and such that given 
any coordinate patch $\U$ on $Y$ diffeomorphic to an open subset 
of $\mathbb{S}^{n-1}$ and function $\varphi \in C^\infty_c(\U)$, 
we have $(1-\omega)\varphi\,u \in H^s(\R^n)$ where $\R_+\times
\mathbb{S}^{n-1}$ is identified with $\R^n \setminus \{0\}$ via polar 
coordinates. By definition, we have 
\[ \mathcal{K}^{0,0}(Y^\wedge)=H_b^0(Y^\wedge)=L_b^2(Y^\wedge). \]
These spaces have been systematically considered by Schulze in his edge 
calculus, see e.g. \cite{Sz89, Sz91}. For $A\in x^{-\mu}\Psi_b^{\mu}(M)$,
the associated model operator $A_\wedge$ extends as a bounded operator 
$A_\wedge: \mathcal{K}^{s,\alpha}(Y^\wedge)\to 
\mathcal{K}^{s-\mu,\alpha-\mu}(Y^\wedge)$ for every $s, \alpha\in\R$.

\begin{definition}\label{ellipticity}
Let $A\in x^{-\mu}\Psi_b^{\mu}(M)$ and let $\Lambda$ be a sector in $\C$ 
containing the origin. The operator family $A-\lambda$ is said to be 
\emph{parameter-elliptic} on $\Lambda$ with respect to $\alpha\in\R$, 
if and only if
\begin{enumerate}[(a)]
\item
$\sym_\mu(A)(\xi)-\lambda$ is invertible for all $\xi\not=0$
and $\lambda\in\Lambda$,
\item
$A_\wedge-\lambda: \mathcal{K}^{s,\alpha}(Y^\wedge)
\to \mathcal{K}^{s-\mu,\alpha-\mu}(Y^\wedge)$ is
invertible for every $\lambda\in\Lambda$ sufficiently large, and for
some $s\in\R$.
\end{enumerate}
These conditions imply that
$\mathrm{spec}_b(A)\cap \{z\in\C \st \Im z=-\alpha\}=\varnothing$.
\end{definition}

It is worth mentioning that in Definition~\ref{ellipticity} (b) the cone 
Sobolev space  $\mathcal{K}^{s,\alpha}(Y^\wedge)$ cannot be replaced by 
the weighted space $x^{\alpha} H_b^{s}(Y^\wedge)$. This is a
consequence of Proposition~\ref{InvNecessary} and the following
elementary observation. 

\begin{lemma}
Let $A_\wedge$ be as in \eqref{ModelOperator} and let $\mu>0$ be such 
that $\mathcal{K}^{\mu,\mu}(Y^\wedge)$ is a proper subspace of 
$x^{\mu}H_b^{\mu}(Y^\wedge)$. If $A_\wedge:x^{\mu}H_b^{\mu}(Y^\wedge) 
\to L_b^{2}(Y^\wedge)$ is invertible, then $A_\wedge: 
\mathcal{K}^{\mu,\mu}(Y^\wedge) \to L_b^{2}(Y^\wedge)$ is not surjective.
\end{lemma}
\begin{proof}
Let $u\in x^{\mu}H_b^{\mu}(Y^\wedge)\setminus\mathcal{K}^{\mu,\mu}(Y^\wedge)$
and assume that $A_\wedge$ is surjective. Then there exists a function 
$v\in \mathcal{K}^{\mu,\mu}$ such that $A_\wedge v=A_\wedge u\in 
L_b^2(Y^\wedge)$. But $v\in x^{\mu}H_b^{\mu}(Y^\wedge)$ and $A_\wedge$ is 
injective, so $u=v\in \mathcal{K}^{\mu,\mu}$ which contradicts the 
assumption on $u$.
\end{proof}

If $\Lambda\subset\C$ is a sector not containing the positive real axis, 
then every operator $A\in x^{-\mu}\Psi_b^{\mu}(M)$ such that 
$A:x^{\alpha}H_b^s(M)\to x^{\alpha-\mu}H_b^{s-\mu}(M)$ is positive and 
selfadjoint, is parameter-elliptic on $\Lambda$ with respect to $\alpha$.
This follows from Proposition~\ref{InvNecessary}.

\begin{example}[cf. Example 3.3 in \cite{GiLo01}]
Let $M$ be a compact $n$-manifold with boundary and let $g$ be a Riemannian 
metric on $M$ which, near the boundary, coincides with the cone metric 
$dx^2 + x^2 g_Y$, where $g_Y$ is a metric on $Y=\partial M$.
The corresponding measure is of the form $x^n \mass$ for a $b$-measure
$\mass$. Let $\Delta_g$ be the Laplace-Beltrami operator associated to
the metric $g$. This operator is symmetric on $L^2(M,x^n \mass) = 
x^{-n/2} L^2_b(M)$ and therefore, the operator
\begin{equation}\label{LaplaceType} 
  A = -x^{n/2-1} \Delta_g\, x^{-n/2+1} + x^{-2} a^2 
      \in x^{-2}\Psi_b^{2}(M) 
\end{equation}
is symmetric on $x^{-1} L^2_b(M)$ for every real number $a$. 
For $u\in C^\infty_c(M)$ supported near the boundary, we have 
\[ Au=x^{-2}\left((xD_x)^2 - \Delta_Y+\tfrac{(n-2)^2}{4}+a^2\right)u, \]
where $\Delta_Y$ is the Laplacian corresponding to $g_Y$. 
For $a>1$ the boundary spectrum of $x^2 A$ does not intersect the 
strip $\{\sigma\in\C \st |\Im\sigma|<1\}$ so that $A$ with domain
$xH_{b}^2(M)$ is positive and selfadjoint on $x^{-1}L_{b}^2(M)$. 
In particular, $A-\lambda$ is parameter-elliptic with respect to $\alpha=1$ 
on any sector $\Lambda\subset\C$ contained in the resolvent set of $A$.
\end{example}

\begin{example}
Let $A$ be the operator \eqref{LaplaceType}.
If $T\in x^{-1}\Psi_b^{1}(M)$ is symmetric on $x^{-1} L^2_b(M)$, then
the operator $A+T$ with domain $xH_{b}^2(M)\embed H_{b}^1(M)$ is also 
positive and selfadjoint, and therefore parameter-elliptic with respect 
to $\alpha=1$. Observe that $T:H_{b}^1(M)\to x^{-1} L^2_b(M)$ is bounded.
\end{example}


We are now ready to prove that our parameter-dependent operators
capture the resolvent of a cone pseudodifferential operators. We
begin by defining certain index sets that appear in Theorem
\ref{thm:res2}. We define
\begin{multline*} 
\widehat{E}^\pm(\alpha) 
= \big\{ (z + r,k) \st r \in \N_0,\, \tau = \mp iz \in \spec_b(A)+i\mu, \\ 
\qquad 1\leq k + 1\leq \sum_{\ell=0}^r \mathrm{ord}(\tau-i \mu \mp i\ell),
\text{ and } \Re z > \pm (\alpha - \mu)\big\},
\end{multline*}
where the order of a pole $\tau \in \spec_b(A)$ of the inverse of the
conormal symbol $\widehat{A}(\tau)$ is denoted by $\mathrm{ord}(\tau)$.
Setting $\check{E}^\pm(\alpha) = \widehat{E}^\pm(\alpha)
\overline{\cup} \widehat{E}^\pm(\alpha)$ and $E(\alpha) = \N\,
\overline{\cup}(\widehat{E}^+(\alpha)+\widehat{E}^-(\alpha))$, we define 
\begin{equation} \label{E(alpha)} 
\E(\alpha) = (\check{E}^+(\alpha),\check{E}^-(\alpha),E(\alpha),\N_0),
\text{ where } \check{E}^\pm(\alpha) = \widehat{E}^\pm(\alpha)
\overline{\cup}\widehat{E}^\pm(\alpha).
\end{equation}

\begin{theorem}\label{thm:res2}
Let $A\in x^{-\mu}\Psi_b^{\mu}(M)$, $\mu>0$, be such that
$A-\lambda$ is parameter-elliptic on $\Lambda$ with respect to some
$\alpha\in\R$. Then for $\lambda \in \Lambda$ sufficiently large,
\[ A -\lambda: x^\alpha H^s_b(M) \to x^{\alpha -\mu}H^{s-\mu}_b(M) \]
is invertible for any $s\in \R$, and
\[ (A - \lambda)^{-1} \in x^{\mu}\Psi^{ - \mu,- \mu,\mu}_c(M;\Lambda) 
   +\ x^{ \mu}\Psi^{-\infty,\mu, \E(\alpha)}_c(M;\Lambda), \] 
where $\E(\alpha)$ is the index family defined in \eqref{E(alpha)}.
Moreover, for $\alpha=\mu=s$ we have
\begin{equation}\label{BoundedResolvent}
 (A -\lambda)^{-1}: L^{2}_b(M)\to x^\mu H^\mu_b(M) 
\end{equation}
is uniformly bounded in $\lambda$.
\end{theorem}
\begin{proof}
Let $(x,y)\in\U = [0,c)_x\times \R^{n-1}$ be local coordinates
near the boundary of $M$ and let $a_\mu(x,y,\xi)$ denote the totally
characteristic principal symbol of $A$. Given $\eps > 0$, let $\chi \in
C^\infty(\R^n)$ with $\chi(\xi)=0$ for $|\xi|<\eps$ and $\chi(\xi)=1$ 
for $|\xi| > 2\eps$. Let
\begin{equation} \label{bmub}
 b_{-\mu}(x,y,\xi,\lambda)
 = \chi(\xi)(a_\mu(x,y,\xi) - x^\mu\lambda)^{-1}
\end{equation}
Observe that $(x',y',z)$, where $z = (\log (x/x'),y-y')$, are
coordinates on $M^2_b$ near $\Delta_b$. Given $\varphi\in C^\infty_c(\U)$ 
and $\psi(z) \in C^\infty_c(\R^n)$ where $\psi(z)=1$ on a neighborhood of 
$z=0$, define the Schwartz kernel of $B(\lambda)$ by
\[ K_{B(\lambda)} =  \varphi(x',y')\, \psi(z) \int e^{iz\cdot \xi}\, 
b_{-\mu}(x',y',\xi,\lambda) \, \dbar\xi \cdot \mass', \] 
where $\mass' = |(dx'/x')dy'|$. Then, by definition, 
$B(\lambda)\in \Psi^{-\mu,-\mu,\mu}_c(M;\Lambda)$ (cf. Example~\ref{ex}). 
Since the principal $b$-symbol of $Ax^\mu$ is also $a_\mu$, and
\[(a_\mu(x,y,\xi) - x^\mu \lambda)\, b_{-\mu}(x,y,\xi,\lambda)\chi(\xi)
  = 1 + (\chi(\xi) - 1), \] 
the composition properties of the $b$-calculus show that
\begin{equation} \label{res1} 
 (A - \lambda)x^\mu B(\lambda) = (Ax^\mu - x^\mu\lambda) B(\lambda)
 = \varphi - S(\lambda) + T, 
\end{equation} 
where $S(\lambda)\in \Psi^{ -1,-\mu,\mu}_c(M;\Lambda)$, and the Schwartz 
kernel of $T$ is given by
\[ K_T = \varphi(x',y')\, \psi(z) \int e^{iz\cdot \xi} \, (\chi(\xi) - 1) 
\, \dbar\xi \cdot \mass'.\] 
Since $\chi(\xi)= 1$ for $|\xi| > 2\eps$, $\chi(\xi) - 1 = 0$ for $|\xi|
> 2 \eps$, which implies that $T$ is a $b$-pseudodifferential operator of 
order $-\infty$ with a symbol supported in $|\xi| < 2 \eps$ and whose 
Schwartz kernel $K_T \to 0$ in the $C^\infty$ topology as a smooth function 
on $M^2_b$. In particular, the mapping properties of $b$-pseudodifferential 
operators \cite{RBM2} imply that the $L^2_b$ norm of $T$ tends to $0$ as 
$\eps \to 0$. If $\U$ is a coordinate patch on the interior of $M$, a 
similar argument shows that given $\varphi \in C^\infty_c(\U)$, there is a 
$B(\lambda)\in\Psi^{-\mu,-\mu,\mu}_c(M;\Lambda)$ such that \eqref{res1} holds.

Let $\{ \U_i \}_{i=1}^N$ be coordinate patches covering $M$ such
that as in \eqref{res1}, there exists a $B_i(\lambda) \in
\Psi^{-\mu,-\mu,\mu}_c(M;\Lambda)$ satisfying $(A - \lambda) x^\mu
B_i(\lambda) = \varphi_i  - S_i(\lambda) + T_i$, where $S_i(\lambda) 
\in \Psi^{-1,-\mu,\mu}_c(M;\Lambda)$, and where $\varphi_i$ is a smooth
function supported in $\U_i$. 
Setting $B_0(\lambda) = \sum_{i=1}^N B_i(\lambda)
\in \Psi^{-\mu,-\mu,\mu}_c(M;\Lambda)$ and assuming that the
$\varphi_i$ form a partition of unity of $M$, we obtain
\begin{equation} \label{res2} 
(A - \lambda)x^\mu B_0(\lambda) = I - S_0(\lambda) + T,
\end{equation} 
where $T \in \Psi^{-\infty}_b(M)$ and $S_0(\lambda)\in 
\Psi^{-1,-\mu,\mu}_c(M;\Lambda)$. Theorem~\ref{thm:comp1} shows that 
$S_0(\lambda)^j \in \Psi^{-j, - j\mu,\mu}_c(M;\Lambda)$ for each $j$. 
Thus we can choose $S_0'(\lambda)\in \Psi^{-1,-\mu,\mu}_c(M;\Lambda)$ 
such that $S_0'(\lambda) \sim \sum_{j = 1}^\infty S_0(\lambda)^j$, 
where the right-hand side is an asymptotic sum. This implies that
\[ (I - S_0(\lambda))(I + S_0'(\lambda)) = I - R_1(\lambda),
   \quad R_1(\lambda) \in \Psi^{-\infty,-\infty,\mu}_c(M;\Lambda). \] 
Multiplying both sides of \eqref{res2} by $I + S_0'(\lambda)$, we obtain
\[ (A - \lambda)x^\mu B_1(\lambda) = I - S_1(\lambda) + T,\] 
where $B_1(\lambda) = B_0(\lambda) + B_0(\lambda) S_0'(\lambda) 
\in \Psi^{-\mu,-\mu,\mu}_c(M;\Lambda)$ by Theorem~\ref{thm:comp1}, and 
$S_1(\lambda) = R_1(\lambda) - TS_0'(\lambda) \in 
\Psi^{-\infty,-\mu,\mu}_c(M;\Lambda)$ by Theorem~\ref{thm:comp2}.

By choosing $\eps> 0$ sufficiently small, we may assume that $I + T$ is 
invertible. The inverse of $I + T$ is of the form $I + T'$ where $T' \in
\Psi^{-\infty,\beta}_b(M)$ for some $\beta > 0$ that depends on the
width of the strip on which the conormal symbol of $T$ is invertible. 
Moreover, since $\|T\|_{L^2_b} \to 0$ as $\eps \to 0$,
the arguments found in \cite[Ch.\ 5]{RBM2} imply that $\beta > 0$
can be choosing arbitrarily large by choosing $\eps > 0$ smaller.
Choose any $N >> 0$ and let $\eps > 0$ be chosen so that $T' \in
\Psi^{-\infty,2N\mu}_b(M)$. Then multiplying both sides of
the previous displayed equation by $I + T'$, we obtain
\begin{equation} \label{res3}
(A - \lambda)x^\mu B_2(\lambda) = I - S_2(\lambda), 
\end{equation} 
where $B_2(\lambda)=B_1(\lambda)+ B_1(\lambda) T' \in 
\Psi^{-\mu ,-\mu,\mu}_c(M;\Lambda) + \Psi^{-\infty,\mu}_{1,N}(M;\Lambda)$ 
by Theorem~\ref{thm:comp3}, and $S_2(\lambda)=S_1(\lambda)+S_1(\lambda)T' 
\in \Psi^{-\infty,\mu}_{1,N}(M;\Lambda)$ by Lemma~\ref{lem-regbounds} and 
Theorem~\ref{thm:comp3}. By Theorem~\ref{thm-compositemN}, 
$S_2(\lambda)^j \in \OpSpaceB{-\infty, \mu}{j,N}$,
which implies that $S_2'(\lambda) = \sum_{j=1}^{N - 1} S_2(\lambda)^j 
\in \OpSpaceB{-\infty,\mu}{1,N}$ satisfies 
\[ (I - S_2(\lambda))(I + S_2'(\lambda)) = I- S_3(\lambda),\quad 
   S_3(\lambda)= S_2(\lambda)^N \in \OpSpaceB{-\infty,\mu}{N,N}.\]
Multiplying both sides of \eqref{res3} by $I + S_2'(\lambda)$, we obtain
\[ (A - \lambda)x^\mu B_3(\lambda) = I - S_3(\lambda), \] 
where $B_3(\lambda)= B_2(\lambda) + B_2(\lambda) S_2'(\lambda)\in 
\OpSpaceA{-\mu ,-\mu,\mu} +\OpSpaceB{-\infty,\mu}{1,N}$ by 
Theorems~\ref{thm:comp3} and \ref{thm-compositemN}, and by 
Lemma~\ref{Lempsi-inftysub} we have $S_3(\lambda)\in 
\OpSpaceB{-\infty, \mu}{N,N}\subset \OpSpaceB{-\infty,\mu, \E}{N}$ 
where $\E$ is the index family on $\T_d$ given by 
$\E=(\varnothing, \varnothing, \N_0,\N_0)$. 

Finally, using the localized inverse $(A_\wedge-\lambda)^{-1}$ (which 
exists by condition (b) in Definition~\ref{ellipticity}) one can modify
the parametrix $x^\mu B_3(\lambda)$ to get a remainder term that decays 
as $1/|\lambda|$. Then, by means of a standard Neumann series argument, 
this new parametrix can be further refined to obtain the exact resolvent. 
The difficulty is to understand the pseudodifferential structure of the 
resolvent which requires understanding the structure of 
$(A_\wedge-\lambda)^{-1}$. This analysis is rather long
but can be done following the same arguments as in the proof of
Theorem 6.1 {in} \cite{LoRes01}. The conclusion is that 
\[ (A - \lambda)^{-1} \in x^{\mu}\OpSpaceA{-\mu,-\mu,\mu} + 
x^{\mu}\OpSpaceB{-\infty,\mu}{1,N} + x^{\mu}\OpSpaceB{-\infty,\mu, \F}{N},\] 
where $\F = (F_{lb},F_{rb},F_{\ff},\N_0)$, with $F_{lb}> \alpha - \mu$, 
$F_{rb}> -(\alpha-\mu)$, and $F_{\ff} > 0$. Now, since 
$\OpSpaceB{-\infty,\mu}{1,\infty} \subset\OpSpaceA{-\infty,-\mu,\mu}$ 
and $\OpSpaceB{-\infty,\mu,\F}{\infty} \subset \OpSpaceA{-\infty,\mu,\F}$, 
and since $N$ can be chosen arbitrarily large, it follows that
\[ (A - \lambda)^{-1} \in x^{\mu}
\OpSpaceA{-\mu,-\mu,\mu} + x^{\mu}\OpSpaceA{-\infty,\mu, \F}, \]
According to \cite[Th. 5]{RBM2} or \cite[Th. 4.4]{Maz91}, we know that for 
fixed $\lambda$, the resolvent $(A - \lambda)^{-1}$ has expansions at 
$\lb$, $\rb$, and $\ff$, with index sets $\check{E}^+(\alpha)$, 
$\check{E}^-(\alpha)$, and $E(\alpha)$, respectively. It follows that 
$\F$  must equal the index set $\E(\alpha)$ given in \eqref{E(alpha)}. 
This proves the first statement of the theorem.

The norm estimate for \eqref{BoundedResolvent} essentially follows 
from corresponding estimates for $x^\mu B_0(\lambda)$ and 
$(A_\wedge-\lambda)^{-1}$.  First of all, observe that $(A-\lambda)^{-1} 
\in \mathcal{L}(L^{2}_b(M),x^\mu H^\mu_b(M))$ is uniformly bounded in
$\lambda$ if and only if $\|(A-\lambda)^{-1}\|_{\mathcal{L}(L^{2}_b(M))}
=\O(|\lambda|^{-1})$ as $|\lambda|\to\infty$.

That $\|x^\mu B_0(\lambda)\|_{\mathcal{L}(L^{2}_b(M))}=\O(|\lambda|^{-1})$ 
as $|\lambda|\to\infty$ is a consequence of the fact that the Schwartz 
kernel of $B_0(\lambda)$ is locally given by the symbol \eqref{bmub}.  
On the other hand, the norm estimate for $(A_\wedge-\lambda)^{-1}$ on 
$L^{2}_b(Y^\wedge)$ is a direct consequence of its $\kappa$-homogeneity 
properties. More precisely, for every $\varrho>0$ we have
\begin{equation*}
 (A_\wedge-\lambda)^{-1}=\varrho^\mu\kappa_\varrho^{-1}
 (A_\wedge-\varrho^\mu\lambda)^{-1}\kappa_\varrho.
\end{equation*}
Setting $\varrho=|\lambda|^{-1/\mu}$ and using that $\kappa_\varrho$ 
is an isometry on $L^{2}_b(Y^\wedge)$, this gives
\begin{equation*}
 \|(A_\wedge-\lambda)^{-1}\|_{\mathcal{L}(L^{2}_b(Y^\wedge))}
 =|\lambda|^{-1}\|(A_\wedge-\tfrac{\lambda}{|\lambda|})^{-1}\|_{
  \mathcal{L}(L^{2}_b(Y^\wedge))}.
\end{equation*}
Hence $\|(A_\wedge-\lambda)^{-1}\|_{\mathcal{L}(L^{2}_b(Y^\wedge))} 
=\O(|\lambda|^{-1})$ as $|\lambda|\to\infty$.
\end{proof}

Composing the resolvent with itself $N$ times, we obtain
\[ (A - \lambda)^{-N} \in x^{N \mu}\OpSpaceA{-N \mu,-N \mu,\mu} 
   +\ x^{N \mu}\OpSpaceA{-\infty,\mu, \E_N(\alpha)} \] 
where the index family $\E_N(\alpha)$ is defined inductively from the 
index family $\E(\alpha)$ using the composition Theorems~\ref{thm:comp3} 
and \ref{thm-full} (with $N =\infty$ there). Now composing 
$(A-\lambda)^{-N}$ with a $b$-pseudodifferential operator $B$ and using 
the fact that each space on the right is closed under such compositions 
by Theorems~\ref{thm:comp2} and \ref{thm:comp3}, we obtain the following
corollary.

\begin{corollary} \label{cor:res2}
Let $A\in x^{-\mu}\Psi_b^{\mu}(M)$, $\mu>0$, be such that
$A-\lambda$ is parameter-elliptic with respect to some $\alpha$ on
$\Lambda$. Then given any $B \in x^{-\nu}\Psi^{\mu'}_b(M)$,
$\nu, \mu' \in \R$, for $\lambda \in \Lambda$ sufficiently
large, we have for any $N \in \N$,
\[ B(A - \lambda)^{-N} \in x^{N \mu-\nu}
  \OpSpaceA{\mu'-N \mu,-N \mu,\mu} +\ x^{N
  \mu-\nu}\OpSpaceA{-\infty,\mu, \E_N(\alpha)}, \] 
where $\E_N(\alpha)$ and $\F_N(\alpha)$ are the same index families.
\end{corollary}

We finish this section showing that the invertibility condition (b) 
in Definition~\ref{ellipticity} is necessary for the resolvent to be 
uniformly bounded, cf.  \cite[Th. 4.1]{GKM2}. Although we do not discuss 
here the condition on $\sym_\mu(A)$, it can be proved (as in the case of 
a regular operator on a smooth compact manifold, cf. \cite{SeeleyBook}) 
that (a) is also a necessary condition. It implies that
$A_\wedge-\lambda:\mathcal{K}^{s,\alpha}(Y^\wedge)\to 
\mathcal{K}^{s-\mu,\alpha-\mu}(Y^\wedge)$ is Fredholm, 
so its image is closed.

\begin{proposition}\label{InvNecessary}
Let $A\in x^{-\mu}\Psi_b^{\mu}(M)$, $\mu>0$, be such that 
$A-\lambda:x^\mu H^\mu_b\to L^{2}_b$
is invertible for all $\lambda\in\Lambda$ with $|\lambda|>R$ for some 
$R>0$. If the resolvent 
\[ (A-\lambda)^{-1}: L^{2}_b(M)\to x^\mu H^\mu_b(M)\] 
is uniformly bounded in $\lambda$, then $A_\wedge-\lambda:
\mathcal{K}^{\mu,\mu}(Y^\wedge)\to L^2_b(Y^\wedge)$ 
is invertible for every $\lambda\in\Lambda\setminus\{0\}$.	
\end{proposition}
\begin{proof}
The assumptions on $A-\lambda$ and $(A-\lambda)^{-1}$ imply that, if
$u\in x^\mu H^\mu_b(M)$, then
\begin{equation}\label{A-InjEstimate}
  \|(A-\lambda)u\|_{0}\ge C\|u\|_{\mu}
\end{equation}
for some constant $C>0$, where $\|\cdot\|_{\nu}$ denotes the norm in 
$x^\nu H^\nu_b(M)$. From this estimate we will derive the injectivity of 
$A_\wedge-\lambda: \mathcal{K}^{\mu,\mu}(Y^\wedge)\to L_b^2(Y^\wedge)$. 

Let $v\in C_c^\infty(\interior{Y}{}^\wedge)$ and pick $\U_Y\subset M$  
such that $\U_Y\cong [0,\eps)\times Y$ for some $\eps>0$. Let $\varrho>0$ 
be small enough so that $\kappa_{\varrho}^{-1}v\in
\mathcal{K}^{\mu,\mu}(Y^\wedge)$ is supported in $[0,\eps)\times Y$, so it 
can be regarded as a function in $x^\mu H^\mu_b(M)$ supported in $\U_Y$. 
Let $\varphi,\psi\in C_c^\infty(\U_Y)$ be such that $\psi=1$ on 
$\textrm{supp}(\kappa_{\varrho}^{-1}v)$ and $\varphi \psi=\psi$. We have
\begin{align*}
 \|(\varrho^\mu\kappa_\varrho \varphi A\psi \kappa_\varrho^{-1}-\lambda)v\|_0
 &= \varrho^\mu\|\kappa_\varrho(\varphi A-\varrho^{-\mu}\lambda)
    \psi\kappa_\varrho^{-1}v\|_0 \\
 &= \varrho^\mu\|(\varphi A-\varrho^{-\mu}\lambda)
    \psi\kappa_\varrho^{-1}v\|_0 \\
 &\ge \varrho^\mu \|(A-\varrho^{-\mu}\lambda)\kappa_\varrho^{-1}v\|_0
  - \varrho^\mu \|(1-\varphi)A\psi\kappa_\varrho^{-1}v\|_0 
\end{align*}
since $\kappa_\varrho$ is an isometry on $L_b^2$ and $\varphi=1-(1-\varphi)$.
Note that $(1-\varphi)A\psi$ is a smoothing operator, so the second norm
on the right-hand side of the inequality is uniformly bounded in $\varrho$.
On the other hand, for $\varrho<1$ we can apply \eqref{A-InjEstimate} and get
\begin{align*} 
 \|(A-\varrho^{-\mu}\lambda)\kappa_\varrho^{-1}v\|_0 
 &\ge C\|\kappa_\varrho^{-1}v\|_\mu = C \varrho^{-\mu}\|v\|_\mu.
\end{align*}
Thus, for $\varrho$ small, 
\[ \|(\varrho^\mu\kappa_\varrho\varphi A\psi\kappa_\varrho^{-1}-\lambda)v\|_0
   \ge C\|v\|_\mu + \O(\varrho^\mu). 
\]
Taking the limit as $\varrho\to 0$, by \eqref{ModelOperator} we get
\begin{equation}\label{CinfEstimate}
  \|(A_\wedge-\lambda)v\|_0\ge C\|v\|_\mu 
\end{equation}
for every $v\in C_c^\infty(\interior{Y}{}^\wedge)$. Since this space is 
dense in $\mathcal{K}^{\mu,\mu}(Y^\wedge)$, \eqref{CinfEstimate} also 
holds for every $v\in \mathcal{K}^{\mu,\mu}(Y^\wedge)$ and we get the 
injectivity of $A_\wedge-\lambda$ on $\mathcal{K}^{\mu,\mu}(Y^\wedge)$.

Finally, note that the invertibility assumption on $A-\lambda$ implies
the invertibility of the formal adjoint $A^\star-\bar\lambda:L_b^2(M)\to
x^{-\mu}H_b^{-\mu}(M)$. By the previous argument, this implies the
injectivity of $A_\wedge^\star-\bar\lambda:L_b^2(Y^\wedge)\to 
\mathcal{K}^{-\mu,-\mu}(Y^\wedge)$, and consequently, the surjectivity of 
$A_\wedge-\lambda:\mathcal{K}^{\mu,\mu}(Y^\wedge)\to L_b^2(Y^\wedge)$. 
\end{proof}

\section{Asymptotic expansions}
\label{AsymptoticExpansion}

To obtain an asymptotic expansion of $B(A-\lambda)^{-N}$, we will use 
the following known lemmas whose proofs can be found in 
\cite[Appendix A]{LoRes201}. 

\begin{lemma} \label{prop:fund2}
Suppose that $u(x,y)$ is a compactly supported on $[0,1)^2$ with
expansions at $x=0$ and $y = 0$ given by index sets (not necessarily 
$C^\infty$ $E_{lb}$ and $E_{rb}$, respectively. Then the 
function $v(x)$ defined by 
\[ v(x)=\int_0^1 u(x/y,y)\frac{dy}{y}= \int_0^1 u(y,x/y)\frac{dy}{y}, \] 
can be expanded at $x=0$ with index set $E_{lb} \overline{\cup} E_{rb}$.
\end{lemma}

This lemma is a special case of the ``Pushforward Theorem'' due to 
Melrose~\cite{MeR92}. As discussed in \cite{GrGr99}, this theorem is  
related to the ``Singular Asymptotics Lemma'' due to Br{\"u}ning and 
Seeley~\cite{BrSe91}.

\begin{lemma} \label{prop:fund1}
Let $f \in C^\infty(\R_+)$ vanish to infinite order as $x\to \infty$ and 
suppose that for some $a \in \C$, we have
\begin{equation} 
  (x\partial_x - a)f(x) = g(x),\label{xpxfg} 
\end{equation}
where $g(x)$ can be expanded at $x=0$ with index set $E$, not
necessarily a $C^\infty$ index set. Then $f$ has an expansion at
$x=0$ with index set $E \overline{\cup} \{a\}$.
\end{lemma}

We are now ready to prove our main result concerning asymptotic
expansions of resolvents of pseudodifferential cone operators.

\begin{theorem}
Let $A\in x^{-\mu}\Psi_b^{\mu}(M)$, $\mu>0$, be such that $A-\lambda$ is 
parameter-elliptic on $\Lambda$ with respect to some $\alpha$. Then, given 
any $B \in x^{-\beta}\Psi^{\mu'}_b(M)$ with $\beta, \mu'\in\R$, for $N$ 
sufficiently large, $B(A-\lambda)^{-N}:x^{\alpha-\mu} H_b^s(M) \to x^{\alpha
-\mu-\beta} H_b^{s-\mu'}(M)$ is trace class for every $s\in\R$, and
\begin{multline}\label{trresolvent}
 \Tr B(A - \lambda)^{-N} \sim_{|\lambda| \to \infty}
 \sum_{k=0}^\infty \Big\{ a_k + b_k\log\lambda +
 c_k (\log\lambda)^2 \Big\}\, \lambda^{(\mu'+n-k)/\mu - N} \\
  +\ \sum_{k=0}^\infty \Big\{ d_k + e_k \log \lambda\Big\}
 \lambda^{(\beta-k)/\mu - N} +\ \sum_{k=0}^\infty f_k \lambda^{-k - N}.
\end{multline}
Moreover, $b_k = 0$ unless $k\in (\N_0+\mu' + n -\beta)\cup
(\mu\N_0 + \mu' + n)$; $c_k = 0$ unless $k\in\mu\N_0 \cap(\N_0-\beta) 
+ \mu'+n$; and $e_k = 0$ unless $k \in \mu\N_0 + \beta$.
\end{theorem}
\begin{proof} 
By Corollary \ref{cor:res2}, for $\lambda \in \Lambda$ sufficiently
large, we can write
\[ B(A - \lambda)^{-N} = F(\lambda) + G(\lambda), \]
where $F \in x^{N \mu-\beta} \OpSpaceA{\mu'-N \mu,-N \mu,\mu}$ and 
$G \in x^{N \mu-\beta} \OpSpaceA{-\infty,\mu, \E_N(\alpha)}$. 
Hence   $F(\lambda)\in x^{N \mu-\beta}\Psi_b^{\mu'-N \mu}(M)$ and 
$G(\lambda)\in x^{N \mu-\beta} \Psi_b^{-\infty,\E_N(\alpha)}(M)$ for every 
$\lambda$. Thus by their mapping properties and Remark~\ref{embeddings}, 
the operators $F(\lambda)$ and $G(\lambda)$ are both trace class if $N$ 
is large enough. We assume $\mu'-N\mu<-n$. The expansion \eqref{trresolvent} 
will be achieved by expanding $\Tr F(\lambda)$ and $\Tr G(\lambda)$.

{\sc Step 1:} 
We begin by showing that, as $|\lambda| \to \infty$ in $\Lambda$, we have
\begin{equation} \label{Gtrace}
  \Tr G(\lambda) \sim \sum_{k=0}^\infty \alpha_k
  \lambda^{(\beta-k)/\mu - N},\quad \alpha_k \in \C.
\end{equation}
If $\Delta \cong M$ is the diagonal in $M^2$, then $\Tr G(\lambda)
= \int_M K_{G(\lambda)}|_{\Delta}$. By the definition of 
$x^{N \mu-\beta} \OpSpaceA{-\infty,\mu,\E_N(\alpha)}$, on the interior 
of $\Delta$, $K_{G(\lambda)}|_{\Delta}$ vanishes to infinite order as
$|\lambda| \to \infty$. Thus we may assume that $K_{G(\lambda)}|_{\Delta}$ 
is supported in a neighborhood $[0,1)_x\times Y$ of $M$ near $Y$. 
Let $r = |\lambda|^{-1/\mu}$ and $\theta = \lambda/|\lambda|$. Then,
integrating out the variables on $Y$, we can write (for $r\le 1$)
\begin{align*}
  \int_M K_{G(\lambda)}|_{\Delta} 
  &= \int_0^{1/r}\! x^{N \mu - \beta} G(r,\theta,x/r)\frac{dx}{x}\\ 
  &= r^{N \mu - \beta} \int_0^1 x^{N\mu-\beta}G(r,\theta,x)\frac{dx}{x}
     \quad (x \mapsto rx),
\end{align*}
where $G(r,\theta,v)$ is a function smooth in $r$ up to $r = 0$, smooth 
in $\theta$, can be expanded at $v = 0$ with index set $E_{N,\ff}(\alpha)
\ge \mu-N\mu$, and vanishes to infinite order as $v\to\infty$. Since 
$G(r,\theta,v)$ is smooth at $r=0$, as $r \to 0^+$ we have 
\[ \Tr G(\lambda)\sim \sum_{k=0}^\infty g_k(\theta)\,r^{N\mu-\beta+k} \]
for some $g_k(\theta)$, smooth in $\theta$. Since $r=|\lambda|^{-1/\mu}$ 
and $G(\lambda)$ is holomorphic in $\lambda$, this expansion is really an 
expansion in $\lambda$ (cf. \cite[Prop.\ 5.1]{LoRes201}), which proves 
\eqref{Gtrace}.

It remains to prove an asymptotic of $\Tr F(\lambda)$ as $|\lambda|\to
\infty$. If $\varphi \in C^\infty(M)$ vanishes near the boundary $Y$, 
then the trace of $\varphi F(\lambda)$ can be analyzed using techniques 
similar to \cite{GrSe95}, for instance. The result is 
\begin{equation*} 
 \Tr \varphi F(\lambda) \sim_{|\lambda| \to \infty}
 \sum_{k=0}^\infty \Big\{a_k+b_k \log \lambda \Big\}\, 
 \lambda^{(\mu'+n-k)/\mu - N} + \sum_{k=0}^\infty f_k \lambda^{-k - N},
\end{equation*}
where $b_k=0$ unless $k\in(\N_0+\mu'+n-\beta)\cup(\mu\N_0+\mu'+n)$. 

Thus it suffices to assume that $F(\lambda)$ is supported in a collar 
$[0,1)_x \times Y$. By taking a partition of unity of $Y$, we may assume 
that $F(\lambda)$ is supported in a coordinate neighborhood in the $Y$ 
factor. Also, as with the expansion for $\Tr G(\lambda)$, we only need to 
prove an expansion of the form \eqref{trresolvent} with $\lambda$ replaced 
by $r^{-\mu}$ and coefficients that depend smoothly on 
$\theta=\lambda/|\lambda|$. In other words, we will prove the expansion 
\begin{multline}\label{Ftrace}
 \Tr F(\lambda) \sim_{r\to 0^+}
 \sum_{k=0}^\infty \Big\{a_k(\theta) + b_k(\theta)\log r +
 c_k(\theta) (\log r)^2 \Big\}\, r^{N\mu-\mu'-n+k} \\
  +\ \sum_{k=0}^\infty \Big\{ d_k(\theta) + e_k(\theta)\log r\Big\}
  r^{N\mu-\beta+k} +\ \sum_{k=0}^\infty f_k(\theta) r^{(k+N)\mu}.
\end{multline}
Note that $\theta$ appears only as a parameter, so we may assume without 
loss of generality that $\Lambda = [0,\infty)$. We will complete our proof 
in two more steps.

{\sc Step 2:} 
We reduce \eqref{Ftrace} to an application of Lemma \ref{prop:fund2}. 
Using the definition of $x^{N \mu - \beta} \OpSpaceA{\mu'-N \mu,-N \mu,\mu}$ 
and integrating out the $Y$ factor of $[0,1) \times Y$, we can write 
\[ \Tr F(\lambda) = \int_0^1 \hspace{-.4em} \int_{\R^n} x^{N \mu - \beta}
   a(x,\xi,x^\mu \lambda)\, \dbar\xi \frac{dx}{x}, \] 
where $a(x,\xi,\lambda) \in C_c^\infty \big([0,1)_x, 
S_{r,c\ell}^{\mu'-N\mu,-N\mu,\mu}(\R^n_\xi;\Lambda)\big)$. 
By assumption, $\mu'- N\mu<-n$, so the integral in $\xi$ is absolutely 
convergent. If $r = \lambda^{-1/\mu}$, then
\[ \Tr F(\lambda) = \int_0^1 A(x,r/x) \frac{dx}{x}, \]
where $A(x,z)= x^{N\mu-\beta} \int_{\R^n} a(x,\xi,z^{-\mu})\,\dbar\xi$. 
Let $\varphi \in C^\infty(\R_+)$ be such that $\varphi(z) = 1$ for 
$z \leq 1$ and $\varphi(z) = 0$ for $z \geq 2$. Then, for $r\le 1$, 
\begin{equation} \label{TrF}
  \Tr F(\lambda)= \int_0^1 \varphi(r/x)\, A(x,r/x)\frac{dx}{x} 
  + \int_0^{1/r}\! (1-\varphi(r/x))\,A(x,r/x) \frac{dx}{x}. 
\end{equation} 
We analyze the asymptotics of each integral as $r\to 0^+$. For the second 
integral, we make the change of variables $x \mapsto rx$, which gives
\begin{equation} \label{int1}
  \int_0^{1/r}\! (1-\varphi(r/x))\, A(x,r/x) \frac{dx}{x} 
 = \int_0^1 (1-\varphi(1/x))\, A(rx,1/x) \frac{dx}{x}. 
\end{equation} 
Since $A(rx,1/x)=(r x)^{N\mu-\beta} \int_{\R^n} a(rx,\xi,x^{\mu})\,
\dbar\xi$ and $N\mu-\beta>0$, the integral \eqref{int1} converges 
absolutely. Moreover, since $a(x,\xi,\lambda)$ is smooth at $x=0$, 
\eqref{int1} has an expansion at $r=0$ with index set $N\mu-\beta+\N_0$. 
Thus the second integral in \eqref{TrF} contributes an expansion of the 
form given by the second sum in \eqref{Ftrace}.

It remains to analyze the asymptotics of the first integral in \eqref{TrF}. 
Note that $A(x,z)$ has an expansion at $x=0$ with index set 
$N\mu-\beta+\N_0$ since $a(x,\xi,\lambda)$ is smooth at $x=0$. Thus, as 
$\varphi(z)\,A(x,z)$ is compactly supported in $z$ and $x$, we can apply 
Lemma~\ref{prop:fund2}: If $A(x,z)$ has an expansion at $z=0$ with some 
index set $E$, then the first integral in \eqref{TrF} has an expansion as 
$r=\lambda^{-1/\mu} \to 0^+$ with index set $E\overline{\cup} (N\mu-\beta
+\N_0)$ (see \eqref{extunion} for the definition of $\overline{\cup}$). 
In the following step we will show that 
\[ E=(\mu N+\mu\N_0)\overline{\cup} (N \mu-\mu'-n +\N_0). \]

\medskip
{\sc Step 3:} 
Since the asymptotics of $A(x,z)$ at $z=0$ do not depend on $x$, we may 
omit the $x$ variable. Thus it suffices to determine the asymptotics of 
\begin{equation*}
  A(z)= \int_{\R^n}a(\xi,z^{-\mu}) \, \dbar\xi \quad\text{at } z=0, 
\end{equation*}
where $a(\xi,\lambda)\in S_{r,c\ell}^{\mu'-N\mu,-N\mu,\mu}(\R^n;\Lambda)$.
Let $\chi(\xi) \in C^\infty(\R^n)$ be such that $\chi(\xi) = 0$ near 
$\xi = 0$ and $\chi(\xi) = 1$ outside a neighborhood of $0$. Then, given 
$L \in \N$, expanding $a(\lambda,\xi)$ in its homogeneous components, we 
can write
\[ A(z) = \sum_{k=0}^{L-1} A_k(z) + R_L(z), \]
where for each $k$, 
\begin{equation}\label{A_k}
 A_k(z)= \int_{\R^n}\chi(\xi)\,a_k(\xi,z^{-\mu})\,\dbar\xi
\end{equation} 
with $a_k(\xi,\lambda)$ anisotropic homogeneous of degree 
$\mu'-N \mu-k$, and where $R_L(z) =\int r_L(\xi,z^{-\mu})\,\dbar\xi$ 
with $r_L(\xi,\lambda)\in S_{r}^{\mu'-N\mu-L,-N\mu,\mu}(\R^n;\Lambda)$.  
In particular, $r_L(\xi,z^{-\mu}) = z^{\mu N}\widetilde{r}_L(\xi,z^\mu)$ 
where $\widetilde{r}_L(\xi,w)$ is smooth at $w=0$ and satisfies estimates 
of the form \eqref{estimater}. These estimates imply that $R_L(z)$ can be 
expanded to higher and higher order at $z=0$ with index set $\mu N+\mu\N_0$ 
as $L$ is chosen larger and larger. Thus it suffices to analyze the 
asymptotics of each $A_k(z)$ at $z=0$.

Recall that $a_k(\xi,\lambda)$ has the following properties:
\begin{itemize}
\item 
$a_k(\delta\xi,\delta^\mu\lambda)= \delta^{\mu'-N\mu-k} a_k(\xi,\lambda)$
for $\delta >0$,
\item 
$a_k(\xi,z^{-\mu}) = z^{\mu N}\widetilde{a}_k(\xi,z^\mu)$ with 
$\widetilde{a}_k(\xi,w)$ smooth at $w=0$. 
\end{itemize}

Now, making the change of variables $\xi \mapsto z^{-1}\xi$ in
\eqref{A_k} and using the homogeneity properties of $a_k$, we get
\[ A_k(z) = z^{N \mu-\mu'-n+k} 
   \int_{\R^n}\chi(\xi/z)\, a_k(\xi,1)\, \dbar\xi.\] 
Let $\gamma=N \mu-\mu'-n+k$.
Since $\big(z\D_z -\gamma\big)z^{\gamma}=0$ and 
$z\D_z \chi(\xi/z) = -(\xi\cdot \D_\xi \chi)(\xi/z)$ where 
$\xi\cdot\D_\xi = \sum \xi_j \D_{\xi_j}$, we have
\begin{align*} 
 \big(z\D_z - \gamma\big) A_k(z)
 &= - z^{\gamma}\hspace{-.3em} \int_{\R^n}(\xi\cdot \D_\xi \chi)(\xi/z)\, 
 a_k(\xi,1) \, \dbar\xi \\ 
 &= -\int_{\R^n}(\xi\cdot \D_\xi \chi)(\xi)\,a_k(\xi,z^{-\mu})\,\dbar\xi \\ 
 &= -z^{\mu N}\int_{\R^n}(\xi\cdot \D_\xi \chi)(\xi)\,
     \widetilde{a}_k(\xi,z^{\mu})\,\dbar\xi 
\end{align*} 
by means of the change $\xi \mapsto z \xi$ and due to the properties of 
$a_k$. Since the function $(\xi\cdot \D_\xi \chi)(\xi)$ is supported in a 
compact subset of $\R^n \setminus \{0\}$, the last integral is absolutely 
convergent and so it can be expanded at $z=0$ with index set 
$\mu N+\mu\N_0$. Hence, Lemma~\ref{prop:fund1} implies that $A_k(z)$ can 
be expanded at $z=0$ with index set $(\mu N + \mu\N_0)\overline{\cup} 
(N \mu-\mu'-n+k)$. Thus, as $A(z)$ is an asymptotic sum of the $A_k$'s, 
$A(z)$ itself can be expanded at $z=0$ with index set $(\mu N+\mu\N_0)
\overline{\cup} (N \mu-\mu'-n + \N_0)$. This completes the proof of
the theorem.
\end{proof}

Let $\Lambda$ be a sector of the form 
\[ \Lambda=\{\lambda \in \C\,|\, \eps_0 \leq \arg(\lambda) 
   \le 2\pi-\eps_0 \;\text{ for some } 0 <\eps_0< \pi/2\}. \] 
Let $A\in x^{-\mu}\Psi_b^{\mu}(M)$, $\mu>0$, 
be such that $A-\lambda$ is parameter-elliptic on $\Lambda$ with respect 
to some $\alpha$. Then the heat operator of $A$ can defined as the 
Cauchy integral
\begin{equation} \label{heatkernel}
e^{-tA} = \frac{i}{2 \pi} \int_\Upsilon e^{-t\lambda} \, (A -
\lambda)^{-1} \, d \lambda,
\end{equation} 
where $\Upsilon$ is an counter-clockwise contour in $\Lambda$ of the form
\[ \Upsilon = a + \{ \lambda \in \C \st \arg(\lambda) = \delta\
  \mbox{or\ }\arg(\lambda) = 2 \pi - \delta\},\ \ a < 0,\
  \eps_0 < \delta< \pi/2. \] 
Integrating by parts $N-1$ times, we can rewrite \eqref{heatkernel} as
\begin{equation} \label{hk2ibp} 
 e^{-tA}= \frac{i}{2\pi}\frac{(-t)^{-N+1}}{(N-1)!} 
 \int_\Upsilon e^{-t\lambda} \, (A -\lambda)^{-N} \, d \lambda. 
\end{equation} 

The asymptotic analysis from \cite[Sec.\ 4.6]{BBlN-HaR86} applied to 
the expansion \eqref{trresolvent} induces the following theorem.

\begin{theorem}
Let $A\in x^{-\mu}\Psi_b^{\mu}(M)$, $\mu>0$, be such that $A-\lambda$  
is parameter-elliptic on $\Lambda$ with respect to some $\alpha$. 
Then, given any $B \in x^{-\beta}\Psi^{\mu'}_b(M)$ with $\beta, \mu'\in\R$, 
the operator $Be^{-tA}$ is trace class for $t > 0$, and
\begin{multline}
\Tr B e^{-tA} \sim_{t \to 0^+}
 \sum_{k=0}^\infty \Big\{ \alpha_k + \beta_k \log t +
 \gamma_k (\log t)^2 \Big\}\,
 t^{(k - \mu' - n)/\mu} \label{trBhk} \\
  +\ \sum_{k=0}^\infty \Big\{ \delta_k + \eps_k \log t \Big\}
 t^{(k- \beta)/\mu} +\ \sum_{k=0}^\infty
 \kappa_k\, t^{k}.
\end{multline}
Moreover, $\beta_k = 0$ unless $k \in (\N_0 + \mu' + n-\beta)
\cup(\mu\N_0 + \mu' + n)$; $\gamma_k = 0$ unless $k \in \mu\N_0
\cap(\N_0-\beta)+\mu'+n$; and $\eps_k = 0$ unless $k\in \mu\N_0+\beta$.
\end{theorem}

Now suppose that $(A-\lambda)^{-1}$ exists on a neighborhood of
$\Lambda$. Then as in \cite{GiLo01} one can show that the complex
power $A^z$ of $A$ exists and defines an entire family of
$b$-pseudodifferential operators satisfying $A^z A^w = A^{z+w}$
for $z,w \in \C$. Using the following formula for the complex
powers in terms of the heat operator
\[
A^z  = \frac{1}{\Gamma(-z)} \int_0^\infty t^{-z} e^{-tA}\,
\frac{dt}{t},\quad \Re z << 0,
\]
we can write 
\[ \zeta_A(z):=\Tr A^z = \frac{1}{\Gamma(-z)}\M(f)(-z), \]
where $\M(f)(z)$ is the Mellin transform of the function
$f(t)= \Tr(e^{-tA})$. Applying the results on the poles of Mellin
transforms found in \cite[Sec.\ 4.3]{BBlN-HaR86}, using the expansion 
\eqref{trBhk} of $\Tr(e^{-tA})$ as $t\to 0$, plus the fact that 
$1/\Gamma(-z)$ vanishes for $z\in \N_0$, we obtain the following theorem.

\begin{theorem}[Analyticity of the Zeta Function] \label{main3} 
The zeta function $\zeta_A(z)$ is holomorphic for $\Re z< -n/\mu$; and 
extends to be meromorphic on the whole complex plane, with (possible) 
simple poles on the set $\big\{ \frac{k-n}{\mu}\st k\in \N_0 \big\}$ 
and with (possible) triple poles on the set $\big\{\frac{k}{\mu}\st 
k\in \N_0,\, \frac{k}{\mu} \not\in\N_0 \big\}$.
\end{theorem}

\section{Properties of the index}
\label{IndexFormula}

At last we consider the problem of finding the index of the closed
extensions of a $b$-elliptic differential cone operator $A$ and give 
a formula for the index of its closure.  To this end we first prove that, 
for the purpose of index calculations, some significant simplifications 
can be made.  In fact, one can reduce the problem to the case where the 
operator has coefficients independent of $x$ near $\partial M$, and even 
more, one can assume $\Dom_{\min}(A)$ to be a weighted Sobolev space. 
These results show that simplifying assumptions made by various authors 
in the past can indeed be used without lost of generality.

\subsection*{Invariance properties of the index}
Let $M$ be a smooth compact manifold with boundary.
Let $A\in x^{-\mu}\Diff_b^m(M)$ be $b$-elliptic, $\mu>0$.
We regard $A$ as an unbounded operator
$A:C_c^\infty(M)\subset x^{\nu}L^2_b(M)\to x^{\nu}L^2_b(M)$
and denote by $\Dom_{\min}(A)$ the domain of the closure of $A$.
It is convenient to assume $\nu=-\mu/2$; we can always reduce to this
case by conjugation with $x^{\nu+\mu/2}$. It is known 
(cf.~\cite{GiMe01,Le97}) that every closed extension $A_\Dom$ of $A$ 
on $x^{-\mu/2}L^2_b(M)$ is Fredholm with index
\begin{equation*}
 \Ind A_\Dom = \Ind A_{\Dom_{\min}} + \dim \Dom/\Dom_{\min}.
\end{equation*}
Note that $\dim \Dom/\Dom_{\min}$ is completely determined by the boundary
spectrum of $A$. In this section we will give an analytic formula for the
index of $A_{\Dom_{\min}}$ using the heat trace asymptotics obtained in
the previous section.
 
We shall need the following lemma which also establishes the notation.

\begin{lemma} \label{NormEquiv} 
On $\Dom_{\min}(A)$, for $\eps>0$ small enough, the operator norm
\begin{equation*}
\|u\|_A = \|u\|_{x^{-\mu/2}L^2_b} + \|Au\|_{x^{-\mu/2}L^2_b}
\end{equation*} and the norm
\begin{equation*}
\|u\|_{A,\eps} = \|u\|_{x^{\mu/2-\eps}L^2_b} +
\|Au\|_{x^{-\mu/2}L^2_b}.
\end{equation*} are equivalent.
\end{lemma}
\begin{proof} 
Recall that the embedding $x^{\mu/2-\eps}L^2_b\embed x^{-\mu/2}L^2_b$ is 
continuous for $\eps<\mu$. The equivalence of the norms follows from the 
continuity of $(\Dom_{\min}(A),\|\cdot\|_A)\embed x^{\mu/2-\eps}L^2_b$ 
which is a consequence of the closed graph theorem.
\end{proof}

Write $D_x=-i\frac{\partial}{\partial x}$. The operator
$A=x^{-\mu}P$ is said to have coefficients independent of $x$ near the
boundary if $(xD_x) P = P (xD_x)$ near $\partial M$.
Write $A=A_0+xA_1$ with $A_0$ having coefficients independent of
$x$ near $\partial M$. Let $\varphi\in C_c^\infty(\R)$,
$\varphi=1$ near $0$. Furthermore, for $\tau>0$ let
$\varphi_\tau=\varphi(x/\tau)$ and let
\[ A_{[\tau]} =\varphi_\tau A_0 + (1-\varphi_\tau)A. \]
Clearly, $A$ and $A_{[\tau]}$ have the same conormal symbol (indicial
family).

\begin{proposition}\label{RedToConst}
For small enough $\tau > 0$ the operator $A_{[\tau]}$ is also $b$-elliptic 
and therefore $\Dom_{\min}(A_{[\tau]})= \Dom_{\min}(A)$. Moreover, as 
$\tau\to 0$, $A_{[\tau]}\to A$ in the graph norm of $A$. Thus, on 
$\Dom_{\min}(A)$,
\[ \Ind A_{[\tau]}=\Ind A \]
for every small $\tau > 0$.
\end{proposition}
\begin{proof} Let $\sym_m(A)$ denote the totally characteristic principal
symbol of $A$. Then,
\begin{align*}
\sym_m(A_{[\tau]})
 &=\varphi_\tau\sym_m(A_0) + (1-\varphi_\tau)\sym_m(A)\\
 &=\varphi_\tau\sym_m(A)+(1-\varphi_\tau)\sym_m(A)-x\varphi_\tau\sym_m(A_1)\\
 &=\sym_m(A)-\tau\tilde\varphi_\tau\sym_m(A_1)
\end{align*} 
with  $\tilde\varphi_\tau=(x/\tau)\varphi(x/\tau)$. Since 
$\tilde\varphi_\tau$ is bounded, $\tau \tilde\varphi_\tau$ is small
for $\tau$ small, and thus the invertibility of $\sym_m(A)$ implies
that of $\sym_m(A)-\tau\tilde\varphi_\tau\sym_m(A_1)$  for such $\tau$.
Hence $A_{[\tau]}$ is $b$-elliptic too. Since $A$ and $A_{[\tau]}$
have the same conormal symbol, we have from \cite[Prop.~4.1]{GiMe01} 
that $\Dom_{\min}(A_{[\tau]})=\Dom_{\min}(A)$.

Further, from the $b$-ellipticity of $A$ it follows that there is
a bounded parametrix $Q:x^\gamma H^s_b\to x^{\gamma+\mu}H^{s+m}_b$
such that
\begin{equation*}
 R=I-QA:x^\gamma H^s_b\to x^{\gamma}H^\infty_b
\end{equation*} is bounded for all $s$ and $\gamma$. Write
\begin{align*}
 A-A_{[\tau]} = x\varphi_\tau A_1
 &= x\varphi_\tau A_1 QA + x\varphi_\tau A_1 R \\
 &= \tau \tilde\varphi_\tau A_1 QA + x\varphi_\tau A_1 R.
\end{align*} Now, $A_1 Q: x^{-\mu/2}L^2_b\to x^{-\mu/2}L^2_b$ is
bounded, so if $u\in \Dom_{\min}(A)$, then
\begin{equation*}
 \|\tau \tilde\varphi_\tau A_1 QAu\|_{x^{-\mu/2}L^2_b} \leq
   c\,\tau\|Au\|_{x^{-\mu/2}L^2_b} \leq c\,\tau\|u\|_A.
\end{equation*} Write $x\varphi_\tau A_1 R= \tau^{1-\eps}(\frac
x{\tau})^{1-\eps} \varphi_\tau\,x^{\eps}A_1 R$ and note that
\begin{equation*}
 x^{\eps}A_1 R: x^{\mu/2-\eps}L^2_b\to x^{-\mu/2}L^2_b
\end{equation*} in continuous. Then using Lemma~\ref{NormEquiv} we get
\begin{equation*}
 \|x \varphi_\tau A_1 Ru\|_{x^{-\mu/2}L^2_b}
 \leq \tilde c\,\tau^{1-\eps}\|u\|_{x^{\mu/2-\eps}L^2_b}
 \leq c\,\tau^{1-\eps}\|u\|_A.
\end{equation*} Altogether,%
\begin{equation*}
 \|(A-A_{[\tau]})u\|_{x^{-\mu/2}L^2_b} \leq C\,\tau^{1-\eps}\|u\|_A
\end{equation*} and thus $A_{[\tau]}\to A$ as $\tau\to 0$.
\end{proof}

\begin{remark}
Norm estimates related to those obtained in the previous proof can be
found in the book by Lesch \cite[Lemma 1.3.10]{Le97}.
\end{remark}

In general, $\Dom_{\min}$ is not a Sobolev space. The problem lies
in the possible presence of elements of $\spec_b(A)$ along the
line $\Im\sigma = -\mu/2$. However, for index purposes, one can
conveniently reduce the analysis to a slightly modified operator
whose closure has a Sobolev space as its domain.

\begin{proposition}\label{RedToSobolev} 
Let $A$ be $b$-elliptic. Let $A_\eps=x^{\eps}A$, and regard it as an 
unbounded operator on $x^{-(\mu-\eps)/2}L^2_b(M)$. If $\eps > 0$ is 
sufficiently small, then 
\[ A_\eps: x^{(\mu-\eps)/2}H^m_b(M)\to x^{-(\mu-\eps)/2}L^2_b(M)\] 
is Fredholm, and
\begin{equation*}
 \Ind A_\eps = \Ind A_{\Dom_{\min}}.
\end{equation*}
\end{proposition}
\begin{proof} 
Write $A=x^{-\mu}P$ with $P\in \Diff^m_b(M)$. Let
$\eta>0$ be so small that there is no $\sigma\in \spec_b(A)$ with
$\mu/2-\eta\le \Im\sigma<\mu/2$ or $-\mu/2<\Im\sigma\le -\mu/2+\eta$. 
The kernel of $A$ on tempered distributions $x^{-\infty}H^{-\infty}_b(M)$ 
is the same as that of $P$, which we denote by $K(P)$. Recall that 
$\Dom_{\max}(A)=\{u\in x^{-\mu/2}L^2_b\st A u\in x^{-\mu/2}L^2_b\}$. 
The kernel $K_{\max}(A)$ of $A:\Dom_{\max}\subset x^{-\mu/2}L^2_b
\to x^{-\mu/2}L^2_b$ consists of those elements of $K(P)$ whose Mellin 
transforms are holomorphic in $\Im\sigma\ge \mu/2$; since these 
elements belong to $x^{-\mu/2}L^2_b$ and $Au\in x^{-\mu/2}L^2_b$. 
That is, their Mellin transforms are holomorphic on $\Im\sigma>\mu/2-\eta$. 
Thus $K_{\max}(A)=K_{\max}(A_\eps)$ if $0<\eps<\eta$.  On the other
hand, the kernel $K_{\min}(A)$ of $A:\Dom_{\min}\subset
x^{-\mu/2}L^2_b\to x^{-\mu/2}L^2_b$ consists of those elements of
$K(P)$ whose Mellin transforms are holomorphic in $\Im\sigma>-\mu/2$; 
indeed in \cite[Prop.~3.6]{GiMe01} it is shown show that 
$\Dom_{\min}=\Dom_{\max}\cap x^{\mu/2-\eta}H^m_b$. Thus if $\eps<\eta$ 
then $K_{\min}(A)=K_{\min}(A_\eps)$. Consequently, 
$\dim K_{\min}(A)=\dim K_{\min}(A_\eps)$.

Finally, note that the formal adjoint of $A$ in $x^{-\mu/2}L^2_b$
is $A^\star = x^{-\mu}P^\star$, where $P^\star$ is the formal
adjoint of $P$ in $L^2_b$, and likewise $A_\eps^\star= 
x^{-\mu+\eps}P^\star$. Now recall that the Hilbert adjoint of
$A_{\Dom_{\min}}$ is $A^\star$ with domain $\Dom_{\max}(A^\star)$,
so the first part of the argument yields $\dim
K_{\max}(A^\star)=\dim K_{\max}(A_\eps^\star)$.
\end{proof}

\subsection*{Index formula}
According to the previous discussion, we can reduce the computation
of the index of the closure of a $b$-elliptic differential operator
$A$ to the case where $A$ has coefficients independent of $x$ near
$\partial M$ and such that
\begin{equation}\label{OperatorA}
 A: x^{\mu/2}H^m_b(M) \to x^{-\mu/2}L^2_b(M)
\end{equation}
is Fredholm. Under these assumptions, we will give a formula for the 
index of $A$ in the spirit of 
\cite{BrSe88,Cho85,FST99,FoxHas,Le97,MeNi96,Pi93,SSS97}
that holds even when $A$ is pseudodifferential.

Recently, Witt \cite{Wi01} proved a factorization theorem for
operator-valued elliptic Mellin symbols.  Using his result, it follows
that there is a cone pseudodifferential operator $B$ with empty boundary
spectrum, and a smoothing Mellin operator $H$, such that $A-B(1 + H)$ is
compact. This implies
\begin{equation*}
 \Ind A = \Ind(B(1+H)) = \Ind B + \Ind(1 + H).
\end{equation*}
Note that $\sym_m(A)=\sym_m(B)$ and $\spec_b(A)=\spec_b(1+H)$. In other
words, this formula separates the index contributions from the totally 
characteristic principal symbol and the boundary spectrum of $A$.

We first discuss the index of $B:x^{\mu/2}H^m_b(M) \to x^{-\mu/2}L^2_b(M)$. 
\begin{lemma}
If $B\in x^{-\mu}\Psi_b^m(M)$ is $b$-elliptic with
$\spec_b(B)=\varnothing$, then
\begin{equation} \label{MSIndex}
\Ind B =\Tr e^{-tB^\star B}-\Tr e^{-tBB^\star} \;\text{ for } t>0, 
\end{equation}
where $B^\star$ is the formal adjoint of $B$.
\end{lemma}
\begin{proof}
In general, the Hilbert space adjoint $B^*$ of $B$ on $x^{-\mu/2}L^2_b$ 
is not equal to but rather a closed extension of the formal adjoint 
$B^\star$. However, since $\spec_b(B)=\varnothing$, we also have 
$\spec_b(B^\star)=\varnothing$ and therefore $\Dom_{\min}(B^\star)=
\Dom_{\max}(B^\star)=x^{\mu/2}H^m_b$. Thus $B^*$ must be equal to 
$B^\star$ with domain $x^{\mu/2}H_b^m$, so \eqref{MSIndex} is nothing 
but the well-known Mckean-Singer identity. 
\end{proof}

The identity \eqref{MSIndex} is not always true because $B^\star$ may be 
different from the Hilbert space adjoint $B^*$. The condition on the 
boundary spectrum of $B$ is what makes it work. The consequence of the 
previous lemma is that since $B^\star$ is a cone pseudodifferential 
operator, we can apply our results from Section~\ref{AsymptoticExpansion} 
to get an asymptotic expansion of the right-hand side of \eqref{MSIndex} 
as $t\to 0$, and obtain 
\begin{equation*}
 \Ind B = \omega(B,B^\star),
\end{equation*}
where $\omega(B,B^\star)$ is the constant term in the expansion.

On the other hand, it follows from Piazza \cite{Pi93} that
\begin{equation*}
 \Ind (1 + H) = - \eta_{\mu/2}(0,1 + \hat{H}),
\end{equation*}
where (cf. also \cite{Me95,MeNi96})
\begin{equation*}
\eta_{\mu/2}(0,1 + \hat{H}) =\frac{1}{2\pi i}\int\limits_{\Im\sigma = -\mu/2}
\Tr\Big(\frac{d}{d\sigma} \hat{H}(\sigma)\,
(1+\hat{H}(\sigma))^{-1}\Big)\, d\sigma.
\end{equation*}

As a consequence, we obtain the following index formula.

\begin{theorem}
Let $A=B(1+H)$ as above. Then the index of \eqref{OperatorA} is given by
\begin{equation*}
\Ind A = \omega(B,B^\star)
-\frac{1}{2\pi i}\int\limits_{\Im\sigma=-\mu/2} \Tr\Big(\frac{d}{d\sigma} 
\hat{H}(\sigma)\, (1+\hat{H}(\sigma))^{-1}\Big)\, d\sigma.
\end{equation*}
\end{theorem}


\end{document}